\documentclass{article}

\usepackage[preprint]{neurips_2022}

\usepackage[utf8]{inputenc} % allow utf-8 input
\usepackage[T1]{fontenc}    % use 8-bit T1 fonts
\usepackage{hyperref}       % hyperlinks
\usepackage{url}            % simple URL typesetting
\usepackage{booktabs}       % professional-quality tables
\usepackage{amsfonts}       % blackboard math symbols
\usepackage{nicefrac}       % compact symbols for 1/2, etc.
\usepackage{microtype}      % microtypography
\usepackage{xcolor}         % colors

\title{Formatting Instructions For NeurIPS 2022}

\author{%
  $^1$Chaoqi Yang, $^2$Cheng Qian, $^1$Navjot Singh, $^3$Cao Xiao, $^{4,5}$M Brandon Westover, \\
  {\bf $^1$Edgar Solomonik, $^1$Jimeng Sun} \\
  $^1$University of Illinois Urbana-Champaign, $^2$IQVIA, $^3$Relativity,\\
  $^4$Massachusetts General Hospital, $^5$Harvard Medical School\\
  $^1$\{chaoqiy2,navjot2,solomon2,jimeng\}@illinois.edu, $^2$alextoqc@gmail.com, \\
  $^3$cao.xiao@relativity.com, $^{4,5}$mwestover@mgh.harvard.edu \\
}

\usepackage{amsmath}
\usepackage{amsthm}
\usepackage{stmaryrd}
\usepackage{mathtools}
\usepackage{stmaryrd}
\usepackage[lined,boxed,commentsnumbered,linesnumbered,ruled]{algorithm2e}
\usepackage{multirow}
\usepackage{booktabs}       % professional-quality tables
\newtheorem{thm}{Theorem}

\DeclareMathOperator*{\argmin}{arg\,min}

\usepackage{color}
\usepackage{xcolor}

\usepackage{enumitem}
\setlist[itemize]{leftmargin=*}
\usepackage{xspace}
\newcommand{\method}{\texttt{ATD}\xspace }
\newcommand{\improve}{\mathrm{impr}}
\newcommand{\init}{\mathrm{init}}
\newcommand{\diag}{\mathrm{diag}}

\newcommand{\D}{\mathbf{D}}

\title{\method: Augmenting CP Tensor Decomposition by \\Self Supervision}

\begin{document}
\maketitle

% \renewcommand{\shortauthors}{xxxx, et al.}

\begin{abstract}
	Tensor decompositions are powerful tools for dimensionality reduction and feature interpretation of multidimensional data such as signals. Existing tensor decomposition objectives (e.g., Frobenius norm) are designed for fitting raw data under statistical assumptions, which may not align with downstream classification tasks. 
	In practice, raw input tensor can contain irrelevant information while data augmentation techniques may be used to smooth out class-irrelevant noise in samples. This paper addresses the above challenges by proposing augmented tensor decomposition (\method), which effectively incorporates data augmentations and self-supervised learning (SSL) to boost downstream classification. To address the non-convexity of the new augmented objective, we develop an iterative method that enables the optimization to follow an alternating least squares (ALS) fashion. We evaluate our proposed \method on multiple  datasets. It can achieve $0.8\%\sim 2.5\%$ accuracy gain over tensor-based baselines. Also, 
	our \method model shows comparable or better performance (e.g., up to $15\%$ in accuracy) over self-supervised and autoencoder baselines while using less than $5\%$ of learnable parameters of these baseline models.
\end{abstract}

\section{Introduction}
Extracting unsupervised features from high-dimensional data is essential in various scenarios, such as physiological signals \citep{cong2015tensor}, hyperspectral images \citep{wang2017hyperspectral} and fMRI \citep{hamdi2018tensor}. 
% Traditional feature reduction techniques, e.g., principle component analysis \citep{eckart1936approximation}, mostly work on matrices. 
Tensor decomposition models are often used for high-order feature extraction \citep{sidiropoulos2017tensor}, among these, CANDECOMP/PARAFAC (CP) decomposition is one of the most popular models. 
The low-rank CP tensor decompositions \citep{kolda2009tensor} assume that the input data is composited by a small set of components, while the reduced features are the coefficients that quantify the importance of each basis, which provides a compact representation.

Under this low-rank assumption, existing tensor decomposition objectives aim to {\em fit individual data samples} with statistical error measures \citep{hong2020generalized,singh2021distributedmemory}, e.g., Frobenius norm or KL-divergence. Though {\em fitness} is an essential principle for feature reduction, common objective functions do not account for downstream tasks, e.g., classification. 

Contrastive self-supervised learning (SSL) \citep{he2020momentum} is recently popular for unsupervised feature learning, which utilizes the class-preserving data augmentations \citep{dao2019kernel} and learns an encoder that can filter out class-irrelevant information. The optimization goal is to {\em enforce alignments} \citep{chen2020simple,wang2020understanding}, ensuring that the anchor sample is closer to the positive sample (which has the same latent class as the anchor sample) than to the negative sample (which is in a different latent class) in the embedding space. In an unsupervised setting, positive samples are given by data augmentations, while negative samples are hard to acquire \citep{arora2019theoretical, he2020momentum,chen2020simple}. Also, previous models are mostly built on deep neural networks, which are often black-box models with tens of thousands of learnable parameters. 

This paper aims to incorporate both the {\em fitness} and {\em alignment} principles into CP tensor decomposition\footnote{Our design may work for other tensor models, such as Tucker decomposition. We leave it to future work.} by augmenting the common fitness objective with a new self-supervised loss. The new self-supervised loss is based on the unbiased estimation of negative samples \cite{chuang2020debiased}, which effectively prevents the sampling bias issue \citep{arora2019theoretical}. The purpose of our design (i.e., introducing SSL into CP tensor decomposition) is to learn class-aware tensor decomposition results for boosting the downstream tensor sample classifications. To address the non-convex subproblems from the new objective, we formulate an iterative method, which solves least squares optimization in each step with closed-form solution.  Main contributions are summarized below.
\begin{itemize}
	\item We propose {\bf augmented tensor decomposition}, named \method, which learns an unsupervised CP structure decomposition by extending the original fitness objective with a self-supervised loss on the contrastiveness of similar and dissimilar tensor samples.
	\item We develop an iterative method to address the non-convex subproblem from the new objective, which enables our algorithm to {\bf follow an alternative least squares fashion}. Our algorithm has {\em asymptotically the same complexity} of each optimization sweep as the common CP-ALS.
	\item We provide {\bf extensive evaluations on four real-world datasets} and compare to recent tensor decomposition models, autoencoder models, and self-supervised models. Our method shows better or comparable prediction performance in various downstream classifications while only requiring much fewer (e.g., less than $5\%$ of) parameters than that of deep learning baselines.
\end{itemize}

\section{Background}
\vspace{-1mm}
{\bf Notation.} 
We use plain letters for scalars, such as $x$ or $X$, boldface uppercase letters for matrices, e.g., $\mathbf{X}$, boldface lowercase letters for vectors, e.g., $\mathbf{x}$, and Euler script letters for tensors, random variables of tensors, and probability distributions, e.g., $\mathcal{X}$. Tensors are multidimensional arrays indexed by three or more indices (modes). For example, an $N$-mode tensor $\mathcal{X}$ is an $N$-dimensional array of size $I_1\times\cdots\times I_N$, where $x_{i_1,\dots,i_N}$ is the element at the $(i_1,\cdots,i_N)$-th position. For matrix $\mathbf{X}$, the $r$-th row and column are $\mathbf{x}^{(r)}$ and $\mathbf{x}_{r}$ respectively, while ${x}_{ij}$ is for the $(i,j)$-th element. $\|\mathbf{X}\|_F$ is the Frobenius norm. For vector $\mathbf{x}$, the $r$-th element is ${x}_r$, and we use $\|\mathbf{x}\|_2$ to denote the vector 2-norm, $\langle\cdot,\cdot\rangle$ for the vector inner product, $\circ$ for the outer product, and $\llbracket\cdot\rrbracket$ for the Kruskal product. Indices in the paper start from $1$, e.g., $\mathbf{x}_{1}$ is the first column of the matrix. 

\subsection{Tensor Modeling}
\vspace{-1mm}
This paper aims to learn tensor bases from unlabeled data and then use the bases to build a feature extractor for downstream classification. Without loss of generality (w.r.t. tensor order), we consider the fourth-order tensor, e.g., a collection of multi-channel Electroencephalography (EEG) signals,
\begin{equation}
	\mathcal{T}=\left[\mathcal{T}^{(1)},\mathcal{T}^{(2)},\dots,\mathcal{T}^{(N)}\right]\in\mathbb{R}^{N\times I\times J\times K}, \notag
\end{equation}
where $\mathcal{T}^{(n)}\in\mathbb{R}^{I\times J\times K}$. The first dimension of $\mathcal{T}$ corresponds to data samples (e.g., one for each patient), while the other three are feature dimensions (e.g., {\em channel by frequency by timestamp}).

\paragraph{Data Model.} To model the above tensor, previous works \citep{kolda2009tensor} assume that 
\begin{itemize}
    \item There are a set of rank-one tensor components $\mathcal{E}=\{\mathcal{E}_1,\dots,$ $\mathcal{E}_R\}$, which are learnable;
    \item The tensor data sample/slice $\mathcal{T}^{(n)}$ admits a low-rank structure and can be represented as a weighted sum of these tensor components $\mathcal{E}$, where $\mathbf{x}^{(n)}$ denotes its coefficient vector;
    \item On top of the low-rank structure, each data sample $\mathcal{T}^{(n)}$ also contains additional element-wise i.i.d. Gaussian noise due to real-world distortion (e.g., physical noise in signal measurements).
\end{itemize}
In the context of downstream classifications, we further assume that each sample $\mathcal{T}^{(n)}$ is semantically associated to one of the latent classes $p\in\{1,\dots,P\}$, and we let $\mathcal{D}_p$ be the sample distribution of class-$p$. Thus, the data sample can be formulated as, $\forall n$,
\begin{equation} \label{eq:classdist1}
	\begin{aligned}
		\mathcal{T}^{(n)} &= \sum_{r=1}^R x^{(n)}_{r}\cdot \mathcal{E}_r+\epsilon^{(n)}~\sim \mathcal{D}_p,~~~p\in\{1,\dots,P\},\\
	\end{aligned}
\end{equation} 
where
\begin{align*}
    \mathcal{E}_r &= \mathbf{a}_{r}\circ\mathbf{b}_{r}\circ\mathbf{c}_{r},~~~r\in\{1,\dots,R\}, \\
		\mathbf{\epsilon}^{(n)}&\sim_\mathrm{i.i.d.}~\mathcal{N}(0,\sigma),~~\textit{where}~\sigma~\textit{is~generally~small}.
\end{align*}
Here $\mathbf{a}_{r},\mathbf{b}_{r},\mathbf{c}_{r}$ are the learnable parameters, which produces the rank-one component $\mathcal{E}_r$, and they are the column vectors of $\{\mathbf{A}\in\mathbb{R}^{I\times R},\mathbf{B}\in\mathbb{R}^{J\times R},\mathbf{C}\in\mathbb{R}^{K\times R}\}$, referred as "bases". 
% In sum, the goal of the paper is to learn tensor bases $\{\mathbf{A},\mathbf{B},\mathbf{C}\}$ given the unlabeled data $\mathcal{T}^{(n)},~\forall n$, so as to extract meaningful coefficient vector $\mathbf{x}^{(new)}$ for new data $\mathcal{T}^{(new)}$.
% , which are parameterized by the bases $\{\mathbf{A}\in\mathbb{R}^{I\times R},\mathbf{B}\in\mathbb{R}^{J\times R},\mathbf{C}\in\mathbb{R}^{K\times R}\}$,
% \begin{equation}
% 	\mathcal{E}_r = \mathbf{a}_{r}\circ\mathbf{b}_{r}\circ\mathbf{c}_{r},~~~r\in[1,\dots,R]. \notag
% \end{equation}
% Here $\mathbf{a}_r,\mathbf{b}_r,\mathbf{c}_r$ are the column vectors and  $\mathcal{E}_r$ is a rank-one tensor.

\vspace{-1mm}
\paragraph{CANDECOMP/PARAFAC Decomposition (CPD).} Given the above tensor model, standard CPD only captures the i.i.d. Gaussian noise by minimizing the negative log-likelihood (NLL), which results in the following standard fitness/reconstruction loss,
\begin{equation}
	\mathcal{L}_{cpd} = \sum_{n=1}^N \left\|{\mathcal{T}}^{(n)}-\llbracket\mathbf{x}^{(n)},\mathbf{A},\mathbf{B},\mathbf{C}\rrbracket\right\|_F^2 =  \left\|\mathcal{T}-\llbracket\mathbf{X},\mathbf{A},\mathbf{B},\mathbf{C}\rrbracket\right\|_F^2. \notag
\end{equation}
Here, the Kruskal product $\llbracket\cdot\rrbracket$ outputs a fourth-order reconstructed tensor from four input factor matrices. For consistency, if the first input is a vector, the output is considered as a third-order tensor. 
% CPD purely seeks a low-rank reconstruction, while there is no mechanism to handle downstream classification. To improve the CPD model, the next section will propose a method to leverage the latent class assignment. Following two {principles of good bases}, our model advances the CPD model by a self-supervised loss, which is designed to: (i) minimize the similarity between feature vectors from different classes; (ii) maximize the similarity of feature vectors within the same class, while no actual labels are needed.

\subsection{Problem Formulation}

CP decomposition seeks a low-rank reconstruction, without special consideration for the downstream task. In this paper, we are motivated to improve the CPD model by exploiting the latent classes (in an unsupervised way) and learn good bases to provide better class-aware features for classification.
%  Our goal is to find good bases, $\mathbf{A}\in\mathbb{R}^{I\times R}$, $\mathbf{B}\in\mathbb{R}^{J\times R}$, $\mathbf{C}\in\mathbb{R}^{K\times R}$, for such tensor data by seeking a rank-$R$ decomposition with {\em coefficient matrix}, $\mathbf{X}\in\mathbb{R}^{N\times R}$, where each data sample is approximated by the linear composition of a set of rank-one components,
% \begin{align}
%     \mathcal{T}^{(n)} \approx \sum_{r=1}^R x^{(n)}_{r}(\mathbf{a}_{r}\circ\mathbf{b}_{r}\circ\mathbf{c}_{r})\in\mathbb{R}^{I\times J\times K}.
% \end{align}
% Here, the components are parameterized by the bases $\{\mathbf{A},\mathbf{B},\mathbf{C}\}$, while the rows of $\mathbf{X}$ are the {\em coefficient vectors} of each sample, quantifying the contribution of the rank-one components.

\vspace{-1mm}
\paragraph{What are Good Bases?} This paper considers two design principles for good bases. The first principle is {\em fitness}, which requires a low-rank tensor reconstruction with the bases.
% , which distills information approximately from raw signals into coefficient vectors. 
Second, data samples associated with the same latent class should be decomposed into similar coefficient vectors, with the bases, while the vectors should be dissimilar if the samples are from different latent classes. This principle is called {\em alignment}, which is important for classification but not considered in the standard tensor decomposition. In this paper, we assess the quality of the learned bases by the performance of downstream classification, where the coefficient vectors (obtained using the bases via decomposition) are the feature inputs (into the downstream classifier).

 \paragraph{Working Pipelines.} To put it succinctly, the paper tackles an unsupervised learning problem while using downstream supervised classification for evaluation. The procedures are briefly outlined: 
\begin{itemize}
	
	\item First, we {\bf learn} the bases $\{\mathbf{A},\mathbf{B},\mathbf{C}\}$ from a large set of unlabeled data. The loss function is developed in consideration of the {\em fitness} and {\em alignment} (defined in the next section) principles.
	\item Then, we {\bf construct} the following feature extractor given $\{\mathbf{A},\mathbf{B},\mathbf{C}\}$. The feature vector of a new sample is obtained by the closed-form solution of the least squares problem ($\alpha>0$ is a hyperparameter),

    \begin{equation} \label{eq:feature_extractor}
	\begin{aligned}
		&\mathbf{f}(\mathcal{T}^{(new)};\mathbf{A},\mathbf{B},\mathbf{C})
		=\argmin_{\mathbf{x}\in\mathbb{R}^{1\times R}}\left(\left\|\mathcal{T}^{(new)}-\llbracket\mathbf{x},\mathbf{A},\mathbf{B},\mathbf{C}\rrbracket\right\|_F^2+\alpha\|\mathbf{x}\|_2^2\right).
	\end{aligned}
	 \end{equation}
	Note that, when $\mathbf{f}(\cdot)$ is applied to a batch of samples, it outputs a coefficient matrix.
	\item Next, we {\bf evaluate} the feature extractor with a set of labeled data. Given $\mathbf{f}(\cdot)$, we first apply it on the labeled data to extract their features and then train an additional logistic regression model (as the downstream classifier) on top of the extracted features, so that the result of classifications will implicitly reflect how good the bases are.
\end{itemize}

\section{Augmented Tensor Decomposition}
\begin{figure*}
	\centering
	\includegraphics[width=0.92\textwidth]{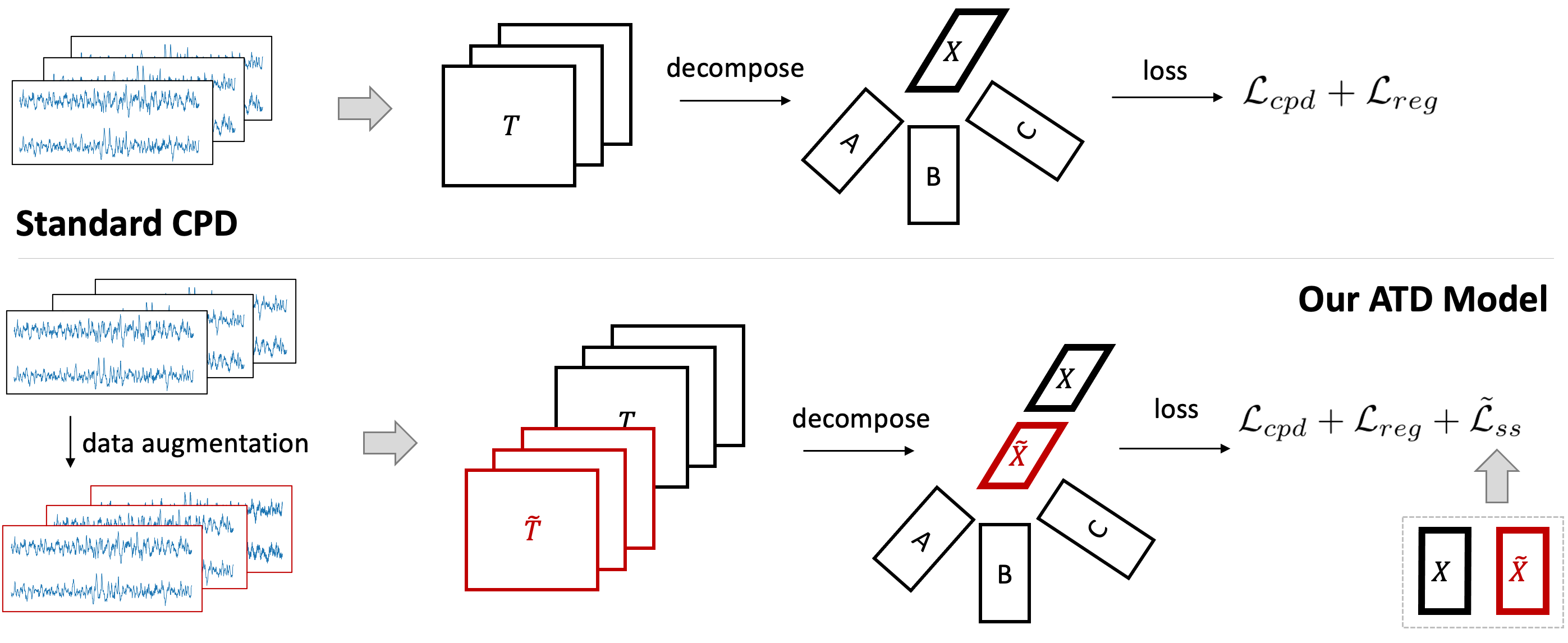}
	\caption{Standard CPD vs Our \method Model}
	\label{fig:framework}
\end{figure*}

We show our model in Figure~\ref{fig:framework}. The design is inspired by the recent popularity of SSL. To exploit the latent class information, we introduce class-preserving data augmentation into CPD model and design self-supervised loss to constrain the learned low-rank features (i.e., the coefficient vectors).  

\noindent{\bf Data Augmentation\footnote{{We provide further discussions and ablation studies on data augmentation in appendix~D.3 and D.5}}.} In general, data augmentation methods are chosen to perturb the raw data while preserving class label information. Given a sample $\mathcal{T}^{(n)}$, we assume that after data augmentation, its perturbation $\tilde{\mathcal{T}}^{(n)}$ preserves the label and admits the component-based representation as in Eqn.~\eqref{eq:classdist1}. 

\subsection{Self-supervised Loss}The design of our self supervised loss corresponds to the {\em alignment} principle, which is based on pairwise feature similarity and dissimilarity. We call a pair of data samples from the same latent class as {\em positive pair}, a pair of samples from different latent classes as {\em negative pair} and a pair of independent samples (two random samples from the dataset) as {\em random pair}. Intuitively, an anchor plus a positive sample composes a positive pair, similarly for negative pairs and random pairs. In this work, our self-supervised loss aims to maximize the feature similarity between positive pairs and minimizes that between negative pairs {\em in an unsupervised way} (no labels during optimization).

Formally, let $\mathcal{X}_p, \mathcal{Y}_p$ be discrete random variables (of tensor samples) distributed as $\mathcal{D}_p,~p\in\{1,\dots,P\}$. We want to minimize the following objective when {\em no class labels are given},
\begin{align*}
    \mathcal{L}_{ss} =&~~\mathcal{L}_{pos} + \lambda\mathcal{L}_{neg},
\end{align*}
where $\lambda\geq 1$ is a hyperparameter and
\begin{equation} \label{eq:ded_pos_neg_loss}
\begin{aligned}
	\mathcal{L}_{pos}&=-\mathbb{E} \left[\mbox{sim}\left(\mathbf{f}\left(\mathcal{X}_p\right),\mathbf{f}\left({\mathcal{Y}_q}\right)\right)\mid p= q\right]
	% =-\mathbb{E}_{i,j\sim d^+} \mbox{sim}\left(\mathbf{f}\left(\mathcal{X}\right),\mathbf{f}\left({\mathcal{Y}}\right)\right)
	,\\
	\mathcal{L}_{neg}&=\mathbb{E} \left[\mbox{sim}\left(\mathbf{f}\left(\mathcal{X}_p\right),\mathbf{f}\left({\mathcal{Y}_q}\right)\right)\mid p\neq q\right].
	% =\mathbb{E}_{i,j\sim d^-} \mbox{sim}\left(\mathbf{f}\left(\mathcal{X}\right),\mathbf{f}\left({\mathcal{Y}}\right)\right)
\end{aligned}
\end{equation}
% \js{I assume the expression after 2nd equal sign is the short form for the first expression? If so, we need to clarify that }
Here $\mathcal{L}_{pos}$ maximizes the feature similarity of positive pairs while $\mathcal{L}_{neg}$ minimizes the feature similarity of negative pairs. $\mathbf{f}(\cdot)$ is the feature extractor, defined in Eqn.~\eqref{eq:feature_extractor}, and the similarity measure is given by cosine distance, parameterized by two random variables,
\begin{equation}
	\mbox{sim}\left(\mathbf{f}\left(\mathcal{X}_p\right),\mathbf{f}\left({\mathcal{Y}_q}\right)\right) = \left\langle\frac{\mathbf{f}(\mathcal{X}_p)}{\left\|\mathbf{f}\left(\mathcal{X}_p\right)\right\|_2},\frac{\mathbf{f}(\mathcal{Y}_q)}{\left\|\mathbf{f}\left(\mathcal{Y}_q\right)\right\|_2}\right\rangle. \notag
\end{equation}

To this end, the key is to implement the self-supervised loss, i.e., Eqn.~\eqref{eq:ded_pos_neg_loss}, in an unsupervised setting. Specifically, we want to construct the sampler of positive pairs and the sampler of negative pairs with unlabeled data only. The sampler of positive pairs (in short, positive sampler) can be easily approximated by data augmentation techniques, which provides "surrogate" positive pairs (given any sample as the anchor, we apply data augmentation methods to generate a perturbed data as positive sample, and then the anchor plus the perturbed data is a positive pair). However, the negative sampler is infeasible, without labels. As a compromise, previous works \cite{he2020momentum,chen2020simple} consider using random sampler to replace the negative sampler given that the random sampler can be easily achieved by picking two independent samples from the dataset. However, this practice is shown to induce sampling bias and hurts the learned representation \cite{arora2019theoretical,chuang2020debiased} since a random pair may be from the same latent class.

\vspace{-2mm}
\paragraph{Construction of Negative Sampler.} In this paper, we consider using the {\em law of total probability} to construct the negative sampler in a statistical way. Formally, assume $r_p$ is the label rate of latent class-$p$ (thus, we have $\sum_{p} r_p=1$), we apply the {law of total probability} and the following holds,
\begin{equation} \label{eq:replacement_trick}
\begin{aligned}
    &\mathbb{E} \left[\mbox{sim}\left(\mathbf{f}\left(\mathcal{X}_p\right),\mathbf{f}\left({\mathcal{Y}_q}\right)\right)\mid p\neq q\right]
    =\sum_{p=1}^{P} r_p\sum_{q\neq p} \frac{r_q}{1-r_p}\mathbb{E}_{p,q} \left[\mbox{sim}\left(\mathbf{f}\left(\mathcal{X}_p\right),\mathbf{f}\left({\mathcal{Y}_q}\right)\right)\right] \\
    =&-\sum_{p=1}^{P} \frac{r_pr_p}{1-r_p}\mathbb{E}_{p,q} \left[\mbox{sim}\left(\mathbf{f}\left(\mathcal{X}_p\right),\mathbf{f}\left({\mathcal{Y}_q}\right)\right)\mid p=q\right]+\sum_{p=1}^{P}\sum_{q=1}^{P} \frac{r_pr_q}{1-r_p}\mathbb{E}_{p,q} \left[\mbox{sim}\left(\mathbf{f}\left(\mathcal{X}_p\right),\mathbf{f}\left({\mathcal{Y}_q}\right)\right)\right] \\
    =&-\mathbb{E} \left[\frac{r_p}{1-r_p}\mbox{sim}\left(\mathbf{f}\left(\mathcal{X}_p\right),\mathbf{f}\left({\mathcal{Y}_q}\right)\right)\mid p=q\right]+\mathbb{E}\left[\frac{1}{1-r_p}\mbox{sim}\left(\mathbf{f}\left(\mathcal{X}_q\right),\mathbf{f}\left({\mathcal{Y}_q}\right)\right)\right].
\end{aligned}
\end{equation}
Here, the usage of $\mathbb{E}[\cdot]$ means that the expectation is taken over four interdependent random variables, i.e., $p,q,\mathcal{X}_p,\mathcal{Y}_q$, while $\mathbb{E}_{p,q}[\cdot]$ means that $p,q$ is fixed and thus it is only taken over two random variables, i.e., $\mathcal{X}_p,\mathcal{Y}_q$. The result shows that the negative sampler can be equivalently replaced by a weighted combination of the random sampler and positive sampler. Here we do not have access to $r_p,~\forall p$ with unlabeled data, this issue is dealt with later. 

% \begin{equation}
%     d^-\left(i,j\right)_{|i} = \frac{1}{1-c(i)} d\left(i,j\right)_{|i} - \frac{c(i)}{1-c(i)} d^+\left(i,j\right)_{|i}. \notag
% \end{equation}
\paragraph{Self-supervised Loss.} Consequently, we can reformulate our self-supervised loss as,
\begin{align}
	\mathcal{L}_{ss} =&~~\mathcal{L}_{pos} + \lambda\mathcal{L}_{neg} \notag\\
	=& -\mathbb{E} \left[\mbox{sim}\left(\mathbf{f}\left(\mathcal{X}_p\right),\mathbf{f}\left({\mathcal{Y}_q}\right)\right)\mid p= q\right] + \lambda\mathbb{E} \left[\mbox{sim}\left(\mathbf{f}\left(\mathcal{X}_p\right),\mathbf{f}\left({\mathcal{Y}_q}\right)\right)\mid p\neq q\right] \label{eq:ss_before_plug}\\
	=&~ \mathbb{E} \left[\frac{\lambda}{1-r_p}\mbox{sim}\left(\mathbf{f}\left(\mathcal{X}_p \right),\mathbf{f}\left({\mathcal{Y}_q}\right)\right)\right] - \mathbb{E}\left[\left(\frac{\lambda r_p}{1-r_p}+1\right)\mbox{sim}\left(\mathbf{f}\left(\mathcal{X}_p \right),\mathbf{f}\left({\mathcal{Y}_q}\right)\right) \mid p = q \right] \label{eq:ss_after_plug}. 
	% {=}&-\mathbb{E}_{i,j\sim d^+} \langle\mathbf{v}(i),\mathbf{v}(j)\rangle +\lambda \mathbb{E}_{i}\left(\mathbb{E}_{j\sim d}\frac{1}{1-c(i)}\langle\mathbf{v}(i),\mathbf{v}(j)\rangle - \mathbb{E}_{j\sim d^+} \frac{c(i)}{1-c(i)}\langle\mathbf{v}(i),\mathbf{v}(j)\rangle\right)\\
\end{align}
From Eqn.~\eqref{eq:ss_before_plug} to Eqn.~\eqref{eq:ss_after_plug}, we use the results in Eqn.~\eqref{eq:replacement_trick}

\paragraph{Two-sided Bound.} The above form still requires label rate information, i.e., $r_p,~\forall p$, we therefore consider using the following approximation to the above loss $\mathcal{L}_{ss}$,
\begin{align}
	\mathcal{L}^{\Theta}_{ss}(\gamma) =&~(\gamma+1)\mathbb{E} \left[\mbox{sim}\left(\mathbf{f}\left(\mathcal{X}_p\right),\mathbf{f}\left({\mathcal{Y}_q}\right)\right)\right] -\mathbb{E} \left[\mbox{sim}\left(\mathbf{f}\left(\mathcal{X}_p\right),\mathbf{f}\left({\mathcal{Y}_q}\right)\right)\mid p= q\right]. \label{eq:L_ss_theta}
\end{align}
Here, $\gamma\geq 0$ is a hyperparameter, while $\mathcal{L}_{ss}$ can be bounded as {(derivations in appendix~E)},
\begin{equation} \label{eq:tightbound}
	C_1\mathcal{L}^{\Theta}_{ss}\left(\frac{\lambda-1}{C_1}\right) \leq \mathcal{L}_{ss} \leq C_2\mathcal{L}^{\Theta}_{ss}\left(\frac{\lambda-1}{C_2}\right),~
C_1=1+\max_{p}\frac{\lambda r_p}{1-r_p},C_2=1+\min_{p}\frac{\lambda r_p}{1-r_p}.
\end{equation}
% \js{can we make some remark about how tight the bound is, and how class balance will affect the bound?}
The equivalence is established when $C_1=C_2$, i.e., the class labels are balanced. To simplify the derivation, we ignore $\lambda$ in the following and let $\gamma$ be a new hyperparameter. Also, the constants $C_1$ and $C_2$ can be absorbed into a weight hyperparameter $\beta$, given in the next section. This bound implies that, an easy-to-compute $\beta\mathcal{L}_{ss}^\Theta(\gamma)$ is often a good approximation of $\mathcal{L}_{ss}$ for some $\beta$. The next section specifies how to compute $\beta\mathcal{L}_{ss}^\Theta(\gamma)$ unsupervisedly as our {\em empirical self-supervised loss}.

\subsection{The Objective of \method} \label{sec:objective}
\paragraph{Empirical Estimator.} We obtain an empirical estimator of $\mathcal{L}^\Theta_{ss}$ with Monte Carlo method. Suppose $\mathcal{T}$ and $\tilde{\mathcal{T}}$ are the input tensor and the augmented tensor respectively, and $\mathbf{X}=\mathbf{f}(\mathcal{T}),\tilde{\mathbf{X}}=\mathbf{f}(\tilde{\mathcal{T}})\in\mathbb{R}^{N\times R}$ are the coefficient/feature matrices. We use the row vectors of $\mathbf{X},\tilde{\mathbf{X}}$ to estimate Eqn.~\eqref{eq:L_ss_theta}. 

The first term $\mathbb{E} \left[\mbox{sim}\left(\mathbf{f}\left(\mathcal{X}_p\right),\mathbf{f}\left({\mathcal{Y}_q}\right)\right)\right]$ essentially requires a random sampler, which is approximated by the average cosine similarity of all possible pairs of non-corresponding row vectors of $\mathbf{X},\tilde{\mathbf{X}}$, while the second term $\mathbb{E} \left[\mbox{sim}\left(\mathbf{f}\left(\mathcal{X}_p\right),\mathbf{f}\left({\mathcal{Y}_q}\right)\right)\mid p= q\right]$ requires a positive sampler, which is estimated by the average cosine similarity of all pairs of corresponding row vectors,
% The random sampler is achieved by picking non-corresponding rows, while the positive sampler is approximated by using the corresponding rows,\js{I am clear what you mean by  non-corresponding  or corresponding rows. Try to rephrase or clarify}
\begin{equation} \label{eq:empirical_ss_loss}
\begin{aligned}
		\tilde{\mathcal{L}}^{\Theta}_{ss}(\gamma) =&~ (\gamma+1)\cdot\frac{1}{N(N-1)}\sum_{n=1}^N\sum_{s\neq n}^N \left\langle\frac{\mathbf{x}^{(n)}}{\|\mathbf{x}^{(n)}\|_2},\frac{\tilde{\mathbf{x}}^{(s)}}{\|\tilde{\mathbf{x}}^{(s)}\|_2}\right\rangle- \frac1N\sum_{n=1}^N \left\langle\frac{\mathbf{x}^{(n)}}{\|\mathbf{x}^{(n)}\|_2},\frac{\tilde{\mathbf{x}}^{(n)}}{\|\tilde{\mathbf{x}}^{(n)}\|_2}\right\rangle\\
		=&~\mbox{Tr}\left(\mathbf{X}^\top \D(\mathbf{X})\mathbf{G}(\gamma)\D(\tilde{\mathbf{X}})\tilde{\mathbf{X}}\right),
\end{aligned}
\end{equation}
where $\D\left(\mathbf{X}\right)= \diag\left(\frac{1}{\|\mathbf{x}^{(1)}\|_2},\cdots,\frac{1}{\|\mathbf{x}^{(N)}\|_2}\right)$ is the row-wise scaling matrix and 
\begin{align}
		\mathbf{G}(\gamma) = \begin{bmatrix}
			-\frac{1}{N} & \frac{\gamma+1}{N(N-1)} & \cdots & \frac{\gamma+1}{N(N-1)} \\
			\frac{\gamma+1}{N(N-1)} & -\frac{1}{N} & \cdots & \frac{\gamma+1}{N(N-1)}\\
			\cdots & \cdots & \cdots & \cdots \\
			\frac{\gamma+1}{N(N-1)} & \frac{\gamma+1}{N(N-1)} & \cdots & -\frac{1}{N}
		\end{bmatrix} \notag.
\end{align}

Note that, the form in Eqn.~\eqref{eq:empirical_ss_loss} is significantly different from tensor discriminant analysis \citep{jia2014low,tao2007general}, which integrates the actual label information as a regularizer and is also different from graph regularized tensor decomposition \citep{cai2010graph}, which incorporates side information, such as geometrical positions \citep{maki2018graph}. Compared to standard noise contrastive estimation (NCE) \citep{gutmann2010noise,chen2020simple} in the area of contrastive SSL, our SSL form in Eqn.~\eqref{eq:empirical_ss_loss} considers a subtraction form instead of the softmax formulation, making it amenable to quadratic optimization, as we will show in Sec.~\ref{sec:stochastic_optimization}.
% , which integrate the given class information as a regularizer to constrain individual factor matrices. Our self-supervised loss is also different from graph regularized tensor decomposition \citep{cai2010graph}, which incorporates domain knowledge and side information, such as geometrical positions \citep{maki2018graph}, and also constrain on individual factors. Though our proposed self-supervised loss can also be viewed as a regularizer, we start from the embedding alignment principle and 

\paragraph{Overall Objective.} According to Eqn.~\eqref{eq:tightbound}, the self supervised loss ${\mathcal{L}}_{ss}$ is bounded by ${\mathcal{L}}^\Theta_{ss}(\gamma)$, while the constants (i.e., $C_1,C_2$) can be absorbed into a weight hyperparameter $\beta$. We let the empirical self-supervised loss, $\tilde{\mathcal{L}}_{ss}{=} \beta\tilde{\mathcal{L}}^\Theta_{ss}(\gamma)$. Our objective follows both the {\em fitness} (i.e., CPD reconstruction loss) and {\em alignment} (i.e., self-supervised loss) principles, while also considering Tikhonov regularization \citep{golub1997tikhonov} to constrain the scale of all parameters,
\begin{equation} \label{eq:objective}
	\mathcal{L} = \mathcal{L}_{cpd} + \mathcal{L}_{reg} + \tilde{\mathcal{L}}_{ss},
\end{equation}
where
\begin{align}
	\mathcal{L}_{cpd} &= \left\|\mathcal{T}-\llbracket\mathbf{X},\mathbf{A},\mathbf{B},\mathbf{C}\rrbracket\right\|_F^2 + \left\|\tilde{\mathcal{T}}-\llbracket\tilde{\mathbf{X}},\mathbf{A},\mathbf{B},\mathbf{C}\rrbracket\right\|_F^2, \notag\\ \mathcal{L}_{reg}&=\alpha\left(\|\mathbf{X}\|_F^2+\|\tilde{\mathbf{X}}\|_F^2+\|\mathbf{A}\|_F^2+\|\mathbf{B}\|_F^2+\|\mathbf{C}\|_F^2\right), \notag\\
	\tilde{\mathcal{L}}_{ss}&= \beta\tilde{\mathcal{L}}^\Theta_{ss}(\gamma) = \beta\mbox{Tr}\left(\mathbf{X}^\top \D(\mathbf{X})\mathbf{G}(\gamma)\D(\tilde{\mathbf{X}})\tilde{\mathbf{X}}\right). \label{eq:objective_ss}
\end{align}
The objective has (i) three hyperparameters, $\gamma,\alpha,\beta>0$; (ii) five factor matrices, $\mathbf{A},\mathbf{B},\mathbf{C}$, $\mathbf{X}$, $\tilde{\mathbf{X}}$. 

\subsection{Alternating Least Squares Optimization} \label{sec:stochastic_optimization}

Ideally, we would like to use alternative least squares (ALS) algorithm \citep{sidiropoulos2017tensor} for optimizing the objective, where each factor matrix is updated in a sequence by solving least squares subproblems. With large scale tensors, we can also resort to mini-batch stochastic ALS \citep{cao2020analysis} to reduce memory footprint of the decomposition. However, the objective in Eqn.~\eqref{eq:objective_ss} is non-convex to $\mathbf{X}$ and $\tilde{\mathbf{X}}$, respectively.  For addressing the non-convex subproblem, we design an iterative method in this section, which only involves solving least square problems.

\paragraph{Addressing Non-convex Subproblem.} We use the subproblem of $\mathbf{X}$ as an example. Given $\mathbf{A}$, $\mathbf{B}$, $\mathbf{C}$, $\tilde{\mathbf{X}}$ fixed, we want to minimize Eqn.~\eqref{eq:objective} by finding the optimal solution, denotd as $\mathbf{X}^*$,
\begin{align}
			\mathbf{X}^{*}\leftarrow& \argmin_{{\mathbf{X}}}\left(\left\|{\mathcal{T}}-\llbracket{\mathbf{X}},\mathbf{A},\mathbf{B},\mathbf{C}\rrbracket\right\|_F^2+\alpha\|{\mathbf{X}}\|_F^2 +\beta \mbox{Tr}\left(\mathbf{X}^\top \D(\mathbf{X})\mathbf{G}(\gamma)\D(\tilde{\mathbf{X}})\tilde{\mathbf{X}}\right)\right). \label{eq:non-convex_prob}
	\end{align}
	
	\begin{itemize}
	    \item First, we are interested to find that the matrix-formed problem in Eqn.~\eqref{eq:non-convex_prob} can be decomposed into row-wise subproblems and to obtain the solution of Eqn.~\eqref{eq:non-convex_prob}, it is suffice to solve each subproblem independently. Let us consider the subproblem of the $k$-th row, which is
	    \begin{align}
		 \argmin_{\mathbf{x}}\left(\left\|{\mathcal{T}}^{(k)}-\llbracket{\mathbf{x}},\mathbf{A},\mathbf{B},\mathbf{C}\rrbracket\right\|_F^2+\alpha\|{\mathbf{x}}\|_F^2 +\beta \mbox{Tr}\left(\frac{\mathbf{x}^\top}{\|\mathbf{x}\|_2}\mathbf{g}^{(k)}\D(\tilde{\mathbf{X}})\tilde{\mathbf{X}}\right)\right), \label{eq:non-convex_prob3}
	    \end{align}
	    where ${\mathcal{T}}^{(k)}$ is the $k$-th slice of ${\mathcal{T}}$, and $\mathbf{g}^{(k)}$ is the $k$-th row of $\mathbf{G}(\gamma)$.
		\item Here, Eqn.~\eqref{eq:non-convex_prob3} is still non-convex. We let the derivative of Eqn.~\eqref{eq:non-convex_prob3} objective to be zero and arrange the terms, which yields,
		\begin{align} \label{eq:non-convex-iter}
		{\mathbf{x}} &= \mathbf{v}_1\mathbf{V}_3 - \frac{\beta\mathbf{v}_2}{2\|\mathbf{x}\|_2}\left(\mathbf{I} - \frac{\mathbf{x}^\top\mathbf{x}}{\|\mathbf{x}\|^2_2}\right)\mathbf{V}_3,
		\end{align}
		where
		\begin{align*}
		\mathbf{v}_1 = \mathbf{T}^{(k)}_1(\mathbf{A}\odot\mathbf{B}\odot\mathbf{C}),~~~~
		\mathbf{v}_2 = \mathbf{g}^{(k)}\D(\tilde{\mathbf{X}})\tilde{\mathbf{X}}, ~~~~
		\mathbf{V}_3 = \left(\mathbf{A}^{\top}\mathbf{A}* \mathbf{B}^{\top}\mathbf{B}*\mathbf{C}^{\top}\mathbf{C}+\alpha\mathbf{I}\right)^{-1}.
		\end{align*}
		Here $\mathbf{x},\mathbf{v}_1,\mathbf{v}_2$ are row vectors and $\mathbf{V}_3$ is a matrix. $\mathbf{T}^{(k)}_1$ is the $1$-mode unfolding of $\mathcal{T}^{(k)}$, $\odot$ is the Khatri-Rao product and $*$ is the Hadamard product (i.e., element-wise product).
		\item We consider the following iterative rule and the fixed point is a solution for Eqn.~\eqref{eq:non-convex-iter}, which is a stationary point of Eqn.~\eqref{eq:non-convex_prob3},
		\begin{align} \label{eq:improved_x}
		    {\mathbf{x}}_{\improve} &= \mathbf{v}_1\mathbf{V}_3 - \frac{\beta\mathbf{v}_2}{2\|\mathbf{x}_{\init}\|_2}\left(\mathbf{I} - \frac{\mathbf{x}^{\top}_{\init}\mathbf{x}_{\init}}{\|\mathbf{x}_{\init}\|^2_2}\right)\mathbf{V}_3.
		\end{align}

		We use an initial guess ${\mathbf{x}}_{\init}$ (obtained by solving Eqn.~\eqref{eq:non-convex_prob3} with while $\beta=0$, which is a least square problem) to start and then we repeat Eqn.~\eqref{eq:improved_x} for each row (i.e., each $k$) and let the improved guess be the initial guess, ${\mathbf{x}}_{\init}\leftarrow {\mathbf{x}}_{\improve}$, to iteratively improve the result.
	\end{itemize}
	
	Theorem~\ref{thm:recursion} (proof in appendix~A) ensures that Eqn.~\eqref{eq:improved_x} converges linearly if $\beta$ is chosen to be sufficiently small. In appendix~B, we verify the liner convergence and also empirically show that one round of Eqn.~\eqref{eq:improved_x} is sufficient in our experiment, where $\beta=2$. The non-convex subproblem of $\tilde{\mathbf{X}}$ can be solved in the same way.

\begin{thm} \label{thm:recursion}
	Given non-zero row vectors $\mathbf{v}_1,\mathbf{v}_2,\mathbf{u}^0\in\mathbb{R}^{d}$, non-singular matrix $\mathbf{V}_3\in\mathbb{R}^{d\times d}$ and $\beta>0$. The sequence $\{\mathbf{u}^t\}$, generated by $\mathbf{u}^{t+1} = \mathbf{v}_1\mathbf{V}_3  - \frac{\beta\mathbf{v}_2}{2\|\mathbf{u}^t\|_2}\left(\mathbf{I}-\frac{\mathbf{u}^{t\top}\mathbf{u}^t}{\|\mathbf{u}^t\|_2^2}\right)\mathbf{V}_3$, satisfies,
	\begin{equation}\small \left\|\mathbf{u}^{t+1}-\mathbf{u}^{*}\right\|_2\leq \frac{\beta(2m+M)\|\mathbf{v}_2\|_2\|\mathbf{V}_3\|_F}{m^3}{\left\|\mathbf{u}^t-\mathbf{u}^{*}\right\|_2}, \notag
	\end{equation}
	where $\mathbf{u}^*$ is the fixed point and $m=\min_t\|\mathbf{u}^t\|_2$, $M=\max_t\|\mathbf{u}^t\|_2$ are the bound of the sequence. With a good $\mathbf{u}^0$ and a sufficiently small $\beta$, we can safely assume $0<m\leq M<\infty$.
\end{thm}

\vspace{-1mm}
\paragraph{Optimization Procedures.} To minimize Eqn.~\eqref{eq:objective_ss}, we alternatively update $\mathbf{A},\mathbf{B},\mathbf{C}$, $\mathbf{X}$, and $\tilde{\mathbf{X}}$, where each subproblem involves only solving least squares problems with closed-form solutions. With large-scale tensors (as in the experiments), we optimizes the factors in mini-batches. Between mini-batches, the basis factors $\mathbf{A, B, C}$ are shared and updated incrementally. We summarize the algortihms in appendix~C. The computation head of the algorithm is matricized tensor times Khatri-Rao product (MTTKRP). The complexity of our optimization algorithm is asymptotically the same as applying CP-ALS, which costs $O(NIJKR)$ to {\em sweep} over the whole tensor once.

\section{Experiments} \label{sec:experiments}
\vspace{-1mm}
This section presents the experimental evaluations. Due to space limitation, additional details, including data augmentations and baseline implementation, are presented in appendix~D.

\vspace{-1mm}
\subsection{Experimental Setup}
\vspace{-1mm}
\paragraph{Data Preparation.} We use four real-world datasets: (i) {\em Sleep-EDF} \citep{kemp2000analysis}, which contains EOG, EMG and EEG Polysomnography recordings; (ii) human activity recognition {\em (HAR)} \citep{anguita2013public} with smartphone accelerometer and gyroscope data; (iii) Physikalisch Technische Bundesanstalt large scale cardiology database {\em (PTB-XL)} \citep{alday2020classification} with 12-lead ECG signals; (iv) Massachusetts General Hospital {\em (MGH)} \citep{biswal2018expert} datasets with multi-channel EEG waves. All datasets are split into three disjoint sets (i.e., unlabeled, training and test) by subjects, while training and test sets have labels. Basic statistics are shown in Table~1.
All models (baselines and our \method) use the same augmentation techniques: (a) jittering, (b) bandpass filtering, and (c) 3D position rotation. 
% \js{where are the detailed definitions of these augmentation. Are all augmentation applied to all methods? Some sounds less intuitive like 3D rotation. We should clarify. }
We provide an ablation study on the augmentation methods in appendix~D.5.
\begin{table*}[t!] \label{tb:datastatistics} \centering
	{\caption{Dataset Statistics}
		\resizebox{0.99\textwidth}{!}{\begin{tabular}{c|cccccccc} 
				\toprule 
				{\bf Name} & {\bf Data Sample Format} & {\bf Augmentations} & {\bf \# Unlabeled ($N$)} & {\bf \# Training} & {\bf \# Test} & {\bf Task} & {\bf \# Class} \\
				\midrule
				{Sleep-EDF} & $I\times J\times K$: 14 $\times$ 129 $\times$ 86 & (a), (b) & 331,208 & 42,803 & 41,078 & Sleep Staging & 5\\
				%   {SHHS} & $I\times J\times K$: ~~ 4 $\times$ 129 $\times$ 55 & 4,082,355 & 235,869 & 217,725 & 260 GB \\
				{HAR} & $I\times J\times K$: 18 $\times$ 33 $\times$ 33 & (a), (b), (c) & 7,352 & 1,473 & 1,474 & Activity Recognition & 6 \\
				{PTB-XL} & $I\times J\times K$: 24 $\times$ 129 $\times$ 75 & (a), (b) & 17,469 & 2,183 & 2,185 & Gender Identification & 2 \\
				{MGH} & $I\times J\times K$: 12 $\times$ 257 $\times$ 43 & (a), (b) & 4,377,170 & 238,312 & 248,041 & Sleep Staging & 5 \\
				\bottomrule
	\end{tabular}}}
	\vspace{-1mm}
\end{table*}
% \js{we should clarify the experiment design. Are unlabeled, training and test disjoint sets? For the patient data, are the split by patients? How are the model build with these 3 datasets? Are the experiment design fair? e.g., are our method using both unlabeled data and labeled training to build the model, while supervised methods only use training data (not unlabeled)? In that case, we are using more data than them so reviewers might argue it is not fair comparison. In any case, we need to clearly describe the evaluation design. What are the classes for each dataset (maybe in appendix?}

\vspace{-1mm}
\paragraph{Baseline Methods.} We include the following comparison models from different perspectives:
\vspace{-1mm}
\begin{itemize}
    % \item {\bf Reference models}: {\em CNN} denotes a convolutional neural network (CNN) based model \citep{biswal2018expert} and {\em $\mbox{CNN}_{Aug}$} denotes the same model with augmented training set. These two models are supervised and work as references. 
    
   \item  {\bf Tensor based models}: ${\bf \method}_{ss-}$ is our variant, which removes the self-supervised loss from the objective in Eqn.~\eqref{eq:objective}; {\em Stochastic alternating least squares (SALS)} applies on the the CPD objective with Tikhonov regularizer, which works on large tensors; {\em Graph regularized SALS (GR-SALS)} augments the objective of SALS with a graph regularizer \citep{maki2018graph,cai2010graph}, define as $\mbox{Tr}\left(\mathbf{X}^\top \mathbf{G}\mathbf{X}\right)$. 

    \item {\bf Self-supervised models}: {\em SimCLR}-$r$ \citep{chen2020simple} and {\em BYOL}-$r$ \citep{grill2020bootstrap} are two popular SSL models with their own objective functions, where $r$ indicates the size of the output representation. 

    \item {\bf Auto-encoder models}: {\em AE}-$r$ denotes a CNN based autoencoder with mean square error (MSE) reconstruction loss, and {\em $\mbox{AE}_{ss}$}-$r$ denotes the same autoencoder model with standard NCE loss in the bottleneck layer, where $r$ is the representation size.
\end{itemize}
% add model variant, $\method_{ss-}$, which removes the self-supervised loss from the objective in Eqn.~\eqref{eq:objective}. We consider the state-of-the-art CPD algorithm with {\em dimension trees} \citep{phan2013fast}, called {\em Fast CPD}. Dimension trees accelerate CP decomposition by reusing intermediate MTTKRP results. However, {\em Fast CPD} algorithm requires to load the full data into memory, and thus it is expensive for large tensors. We also consider the stochastic alternating least squares ({\em SALS}) \citep{maehara2016expected} for the CPD objective, which works on large tensors. 

% In addition, we consider the following deep learning models: (i) two supervised models: a convolutional neural network (CNN) model, {\em Supervised} \citep{biswal2018expert}, and the same model with augmented training set, {\em $\mbox{Supervised}_{Aug}$}; (ii) two self-supervised models: {\em SimCLR}-$r$ \citep{chen2020simple}, {\em BYOL}-$r$ \citep{grill2020bootstrap}; (iii) two autoencoder models: CNN based autoencoder, {\em AE}-$r$, and autoencoder model with self-supervised loss in bottleneck layer, {\em $\mbox{AE}_{ss}$}-$r$, where $r$ is the size of the representation. 
 All models use the unlabeled set to train a feature encoder and use training and test sets to evaluate. Note that, deep neural network models use the same CNN backbone. {In appendix~D.8, we have also compared with two recent supervised tensor learning models, which shows the usefulness of our \method and the large unlabeled set, especially in low-label rate scenarios.}

\vspace{-2mm}
\paragraph{Evaluation and Environments.} We evaluate model performance mainly based on {\em classification accuracy}, where we train an additional logistic classifier \citep{he2020momentum} on top of the feature encoder. Also, for different models, we compare their {\em number of learnable parameters}. The experiments are implemented by {\em Python 3.8.5, Torch 1.8.0+cu111} on a Linux workstation with 256 GB memory, 32 core CPUs (3.70 GHz, 128 MB cache), two RTX 3090 GPUs (24 GB memory each). All training is performed on the GPU. For tensor based models, we use $R=32$ and implement the pipeline in CUDA manually, instead of using {\em torch-autograd}.

\vspace{-1mm}
\subsection{Experimental Results}
\vspace{-1mm}
This section shows the experimental results on downstream classification. We use all the unlabeled data to train the encoder or feature extractor, and use training data (since Sleep-EDF and MGH datasets have enough training samples, we randomly selected a subset of them) for learning a downstream classifier and use all test data. Each experiment is conducted with five different random seeds and the mean and standard deviations are reported. The metrics are the {\em accuracy} and the {\em number of learnable parameters}. All models have 32-dim features in the end, except that for two self-supervised baselines and autoencoder models, which have 128-dim options.

\begin{table*}[t!] \centering
	\caption{Result of Downstream Classification (\%). The table shows that our \method can provide comparable or better performance over all baselines with fewer parameters, especially deep learning models. It also shows the usefulness of considering both {\em fitness} and {\em alignment} as part of the objective.}
	\resizebox{1.01\textwidth}{!}{
		\begin{tabular}{lcccccccc} \toprule 
			\multirow{2}{*}{} & \multicolumn{2}{c}{\bf Sleep-EDF (5,000)} & \multicolumn{2}{c}{\bf HAR (1,473)} & \multicolumn{2}{c}{\bf PTB-XL (2,183)} & \multicolumn{2}{c}{\bf MGH (5,000)}\\
			\cmidrule(r){2-3}  \cmidrule(r){4-5} \cmidrule(r){6-7} \cmidrule(r){8-9}
			& {Accuracy} & {\# of Params.} & {Accuracy} & {\# of Params.} & {Accuracy} & {\# of Params.} & {Accuracy} & {\# of Params.}\\
			\cmidrule{1-9} 
			{\bf Self-sup models:} \\
			%   MoCo & 64.65 $\pm$ 0.615 & 74.99 $\pm$ 0.713 & 74,626 & 65.78 $\pm$ 1.914 & 70.51 $\pm$ 0.647 & 172,104\\
			SimCLR-32 & 84.98 $\pm$ 0.358 & 210,384 & 74.75 $\pm$ 0.723 & 53,286 & 69.25 $\pm$ 0.355 & 200,960 & 67.34 $\pm$ 0.970 & 212,624\\
			SimCLR-128 & {\bf 85.19 $\pm$ 0.358} & 222,768 & 76.69 $\pm$ 0.697 & 65,670 & 68.19 $\pm$ 0.793 & 237,920 & 66.98 $\pm$ 1.331 & 246,608\\
			BYOL-32 & 84.29 $\pm$ 0.405 & 211,440 & 73.71 $\pm$ 2.832 & 54,342 & 65.08 $\pm$ 1.535 & 202,016  & 68.83 $\pm$ 1.168 & 214,736\\
			BYOL-128 & 83.26 $\pm$ 0.337 & 239,280 & 71.79 $\pm$ 1.866 & 82,182 & 65.49 $\pm$ 0.612 & 254,432 & 68.55 $\pm$ 1.339 & 279,632\\
			\midrule
			{\bf Auto-encoders:} \\
			AE-32 & 74.78 $\pm$ 0.723 & 217,216 & 63.13 $\pm$ 0.775 & 62,940 & 59.01 $\pm$ 0.896  & 224,528  & 68.58 $\pm$ 0.427 & 220,088\\
			AE-128 & 75.17 $\pm$ 0.897 & 241,888 & 60.52 $\pm$ 1.604 & 87,612 & 58.29 $\pm$ 0.412 & 298,352 & 67.05 $\pm$ 1.375 & 257,048\\
			%   $\mbox{AE}_{ss}$-32 & 82.20 $\pm$ 0.736 & 47,904 & 68.31 $\pm$ 0.598 & 56,540 & 67.40 $\pm$ 0.899 &  111,188\\
			%   $\mbox{AE}_{ss}$-128 & 82.22 $\pm$ 1.326 & 72,576 & 67.73 $\pm$ 0.714 & 105,788 & 68.43 $\pm$ 1.286 & 148,148\\
			$\mbox{AE}_{ss}$-32 & 80.92 $\pm$ 0.345 & 217,216  & 71.70 $\pm$ 2.135 & 62,940 & 68.47 $\pm$ 0.231 & 224,528 & 71.46 $\pm$ 0.386 &  220,088\\
			$\mbox{AE}_{ss}$-128 & 81.84 $\pm$ 0.259 & 241,888  & 72.43 $\pm$ 1.370 & 87,612 & 68.88 $\pm$ 0.604 & 298,352 & 70.19 $\pm$ 0.617 & 257,048\\
			\midrule
			{\bf Tensor models:} \\
			SALS & 84.27 $\pm$ 0.481 & 7,328 & 91.86 $\pm$ 0.295 & 2,688 & 69.15 $\pm$ 0.483 & 7,296  & 71.93 $\pm$ 0.379 & 9,984\\
			GR-SALS & 84.33 $\pm$ 0.356 & 7,328 &  92.33 $\pm$ 0.282  & 2,688 & 69.02 $\pm$ 0.477 & 7,296  & 72.35 $\pm$ 0.228 & 9,984\\
			$\method_{ss-}$ & 84.19 $\pm$ 0.221 & 7,328 & 92.41 $\pm$ 0.391 & 2,688 & 69.38 $\pm$ 0.612 & 7,296  & 72.78 $\pm$ 0.522 & 9,984\\
			%   AugS-CPD & 86.66 $\pm$ 0.146 & 5,632 & 75.29 $\pm$ 0.282 & 6,016  &  73.13 $\pm$ 0.419 & 9,728\\
			\method & {85.01 $\pm$ 0.224} & 7,328 & {\bf 93.35 $\pm$ 0.357} & 2,688 & {\bf 70.26 $\pm$ 0.523}  & 7,296 & {\bf 74.15 $\pm$ 0.431} & 9,984\\
			\bottomrule
			\multicolumn{9}{l}{*Parenthesis shows the number of training samples. {Our improvements are statistically significant with $p<0.05$ (details in appendix~D.7).}}
	\end{tabular}} \label{tb:comparison_result}
	\vspace{-2mm}
\end{table*}

\subsubsection{Better Classification Accuracy with Fewer Parameters} 
\vspace{-1mm}
From Table~\ref{tb:comparison_result}, \method shows comparable or better performance over the baselines. {We have also reported the running time per epoch/sweep in appendix~D.7 for all models}. Compared to the variant $\method_{ss-}$, our \method can improve the accuracy by $1.0\%\sim 1.9\%$, which shows the benefit of the inclusion of self-supervised loss. SALS and $\method_{ss-}$ have similar performance, while their objectives differ in that $\method_{ss-}$ considers the Frobenius norm of the augmented data. Thus, their accuracy gap is caused by the use of data augmentation. Also, the experiments show that the {\em fitness} and {\em alignment} principles are both important. We observe that with a self-supervised loss (i.e., {\em alignment}), $\mbox{AE}_{ss}$ can give significant improvement over $\mbox{AE}$, while \method shows $\sim 8\%$ accuracy gain over the self-supervised models on MGH dataset, since we can better preserve the data with a reconstruction loss (i.e., {\em fitness}).

Moreover, the table shows that tensor based models require fewer parameters, i.e., less than $5\%$ of parameters compared to deep learning models. On HAR, the deep unsupervised models show poor performance due to (i) they may not optimize a large number of parameters on middle-scale dataset; (ii) movement signals in HAR might have few degrees of freedom, which matches well with the low-rank assumption of tensor methods. On large-scale Sleep-EDF, self-supervised models outperforms \method marginally since they have more parameters thus can capture more information.

\subsubsection{Better Performance in Low-label Rate Scenarios} 
\begin{figure}[t!]
	\centering
	\includegraphics[width=4.75in]{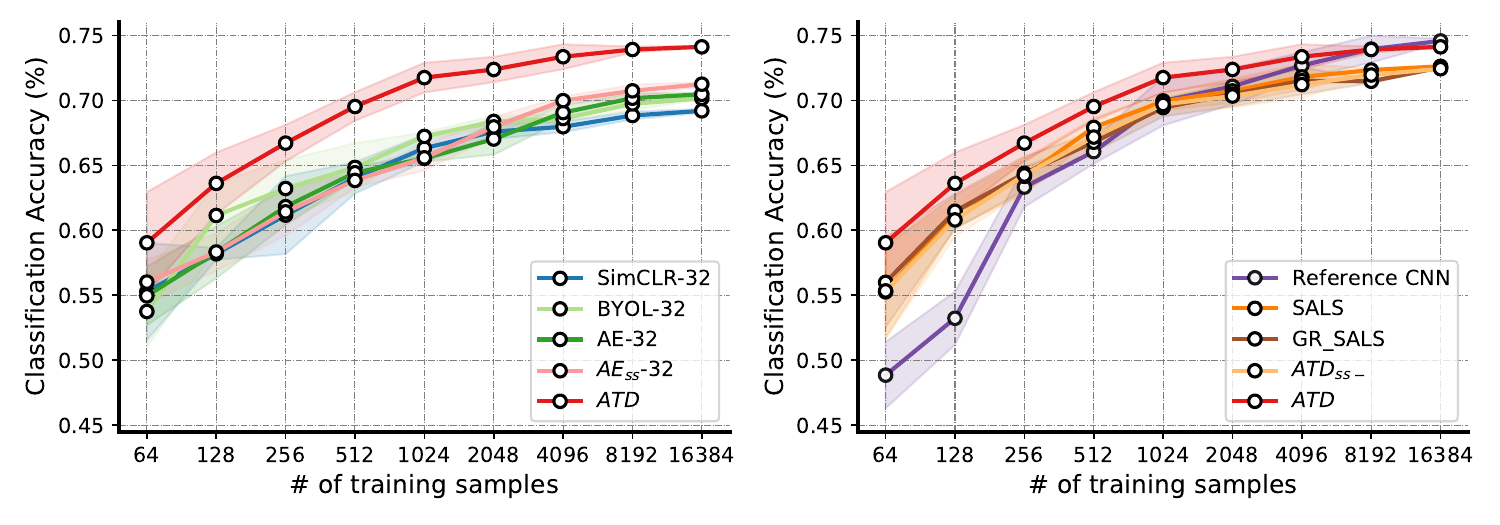}
	\vspace{-3mm}
	\caption{Varying the \# of Training Data} \label{fig:number_of_train_mgh}
\end{figure}

\vspace{-1mm}
On the MGH dataset, we also show the effect of varying the amount of training data in Figure~\ref{fig:number_of_train_mgh}.
We include an end-to-end convolutional neural network (CNN) model based on \citep{biswal2018expert}, called {\em Reference CNN}, which is a supervised model and only uses the training and test sets. To be more readable, we separate the comparison figure into two sub-figures: the left compares our \method model with self-supervised and auto-encoder baselines and the right one compares \method with tensor baselines and the reference model, and the scale of y-axis on two sub-figures are the same.

We find that all unsupervised models outperform the supervised reference CNN model in scenarios with fewer training samples. With more training data, the performance of all models get improved, especially the reference CNN model, which can optimize the encoder and predictive layers in an end-to-end way and finally outperforms our \method when more training samples is available.

\subsubsection{Stable Results with Hyperparameter Variation}
\vspace{-1mm}
For a comprehensive evaluation, we also conduct ablation studies on the effect of the data augmentation methods and on hyperparameters. Due to space limitation, we move the experimental settings and results to appendix~D.6, while summarizing the general conclusions here: (i) with more diverse data augmentation
methods, the final results are relatively better; (ii) with a larger rank $R$, the performance will be better
generally; (iii) our \method is not sensitive hyperparameters $\alpha,\beta,\gamma$.

\section{Conclusion} \label{sec:conclusion}
\vspace{-2mm}
This paper introduces the concept of self-supervised learning for tensors and proposes {\em Augmented Tensor Decomposition} (\method) and the ALS-based optimization. We show that by explicitly contrasting positive and negative samples, the decomposition results are more aligned with downstream classification. On four real-world datasets, we show the advantages of our model over various unsupervised models and with fewer training data, our model even outperforms the reference supervised models.

Compared to deep learning methods, tensor based models are not as flexible in processing multimodal and diverse inputs, such as natural images. However, applying tensor decomposition on the outputs of earlier layers of pre-trained deep neural networks may be a feasible way to address the weaknesses. This direction would be interesting for future work.
% Finally, a case study discussion for how to select good data augmentations is provided. 

\bibliographystyle{apalike}
\bibliography{sample-base}

\clearpage
\appendix
\section{Proof of theorem}
\begin{proof} We provide the proof below (first by submultiplicativity),
{\small \begin{align*}
    \left\|\mathbf{u}^{t+1} - \mathbf{u}^*\right\|_2 &= \left\|\frac{\beta\mathbf{v}_2}{2\|\mathbf{u}^*\|_2}\left(\mathbf{I}-\frac{\mathbf{u}^{*\top}\mathbf{u}^*}{\|\mathbf{u}^*\|_2^2}\right)\mathbf{V}_3 - \frac{\beta\mathbf{v}_2}{2\|\mathbf{u}^t\|_2}\left(\mathbf{I}-\frac{\mathbf{u}^{t\top}\mathbf{u}^t}{\|\mathbf{u}^t\|_2^2}\right)\mathbf{V}_3\right\|_2 \\
    & \leq \left\|\frac{\beta\mathbf{v}_2}{2}\right\|_2\|\mathbf{V}_3\|_F\left\|\frac{1}{\|\mathbf{u}^*\|_2}\left(\mathbf{I}-\frac{\mathbf{u}^{*\top}\mathbf{u}^*}{\|\mathbf{u}^*\|_2^2}\right)-\frac{1}{\|\mathbf{u}^t\|_2}\left(\mathbf{I}-\frac{\mathbf{u}^{t\top}\mathbf{u}^t}{\|\mathbf{u}^t\|_2^2}\right)\right\|_2 \\
    & = \left\|\frac{\beta\mathbf{v}_2}{2}\right\|_2\|\mathbf{V}_3\|_F\left\|\left(\frac{1}{\|\mathbf{u}^*\|_2}-\frac{1}{\|\mathbf{u}^t\|_2}\right)+\left(\frac{\mathbf{u}^{t\top}\mathbf{u}^t}{\|\mathbf{u}^t\|_2^3}-\frac{\mathbf{u}^{*\top}\mathbf{u}^*}{\|\mathbf{u}^*\|_2^3}\right)\right\|_2.
\end{align*}}
Then, we decompose the last term into two parts and have
{\footnotesize \begin{align*}
    \frac{1}{\|\mathbf{u}^*\|_2}-\frac{1}{\|\mathbf{u}^t\|_2} &= \frac{\|\mathbf{u}^t\|_2-\|\mathbf{u}^*\|_2}{\|\mathbf{u}^*\|_2\|\mathbf{u}^t\|_2} \leq \frac{\|\mathbf{u}^t-\mathbf{u}^*\|_2}{\|\mathbf{u}^*\|_2\|\mathbf{u}^t\|_2} \leq \frac{\|\mathbf{u}^t-\mathbf{u}^*\|_2}{m^2}, \\
    \frac{\mathbf{u}^{t\top}\mathbf{u}^t}{\|\mathbf{u}^t\|_2^3}-\frac{\mathbf{u}^{*\top}\mathbf{u}^*}{\|\mathbf{u}^*\|_2^3} & = \left(\frac{\mathbf{u}^{t\top}\mathbf{u}^t}{\|\mathbf{u}^t\|_2^3}-\frac{\mathbf{u}^{*\top}\mathbf{u}^t}{\|\mathbf{u}^t\|_2^3}\right) + \left(\frac{\mathbf{u}^{*\top}\mathbf{u}^t}{\|\mathbf{u}^t\|_2^3} - \frac{\mathbf{u}^{*\top}\mathbf{u}^*}{\|\mathbf{u}^t\|_2^3}\right) + \left(\frac{\mathbf{u}^{*\top}\mathbf{u}^*}{\|\mathbf{u}^t\|_2^3} - \frac{\mathbf{u}^{*\top}\mathbf{u}^*}{\|\mathbf{u}^*\|_2^3}\right) \\
    &= \frac{(\mathbf{u}^{t\top}-\mathbf{u}^{*\top})\mathbf{u}^t}{\|\mathbf{u}^t\|_2^3} + \frac{\mathbf{u}^{*\top}(\mathbf{u}^{t}-\mathbf{u}^*)}{\|\mathbf{u}^t\|_2^3} + \frac{\|\mathbf{u}^*\|^2_2(\|\mathbf{u}^{*}\|_2^3-\|\mathbf{u}^t\|^3_2)}{\|\mathbf{u}^t\|_2^3\|\mathbf{u}^*\|_2^3} \\
    &\leq \frac{\|\mathbf{u}^t-\mathbf{u}^*\|_2}{m^2} + \frac{M\|\mathbf{u}^t-\mathbf{u}^*\|_2}{m^3} + \frac{(M+2m)\|\mathbf{u}^t-\mathbf{u}^*\|_2}{m^3}.
\end{align*}}
The proof is complete by triangular inequality.
\end{proof}

\section{Analysis of the Iterative Rule} In this section, we study how many rounds of the iterative rules are needed to achieve a good classification result. {This study is conducted on HAR and Sleep-EDF with five random seeds.}

% {We first clarify that the "iterative rule" is different from "iteration" in our paper. An "iteration" (e.g., the $l$-th iteration) means five optimization steps (cold start, auxiliary step, three main steps) for one data batch, while one round of the "iterative rule" means applying  Equation~(14) in the main paper once, which is part of the auxiliary step.} 

{\bf Linear Convergence Speed.} First, we study the convergence speed (when $\beta=2$). We consider two scenarios: (Scenario 1) when the bases $\{\mathbf{A},\mathbf{B},\mathbf{C}\}$ are initialized as random matrices; (Scenario 2) when the based $\{\mathbf{A},\mathbf{B},\mathbf{C}\}$ are already learned. We use the average relative difference (of the F-norm) as the convergence measure, i.e., $\frac{1}{N}\sum_{k=1}^N\frac{\|\mathbf{x}^{(k)}_{t+1}-\mathbf{x}^{(k)}_{t}\|}{\|\mathbf{x}^{(k)}_{t}\|}=\frac{1}{N}\sum_{k=1}^N\frac{\|\mathbf{x}^{(k)}_{\improve}-\mathbf{x}^{(k)}_{\init}\|}{\|\mathbf{x}^{(k)}_{\init}\|}$, where $\mathbf{x}^{(k)}$ means the $k$-th row of $\mathbf{X}$ and $t$ means the number of rounds of the iterative rule. We test on $t=1,2,3,4,8$. The comparison is shown in Figure~3. Both scenarios verify the linear convergence speed. 

\begin{figure}[h!] 
    \centering
    \includegraphics[width=0.99\textwidth]{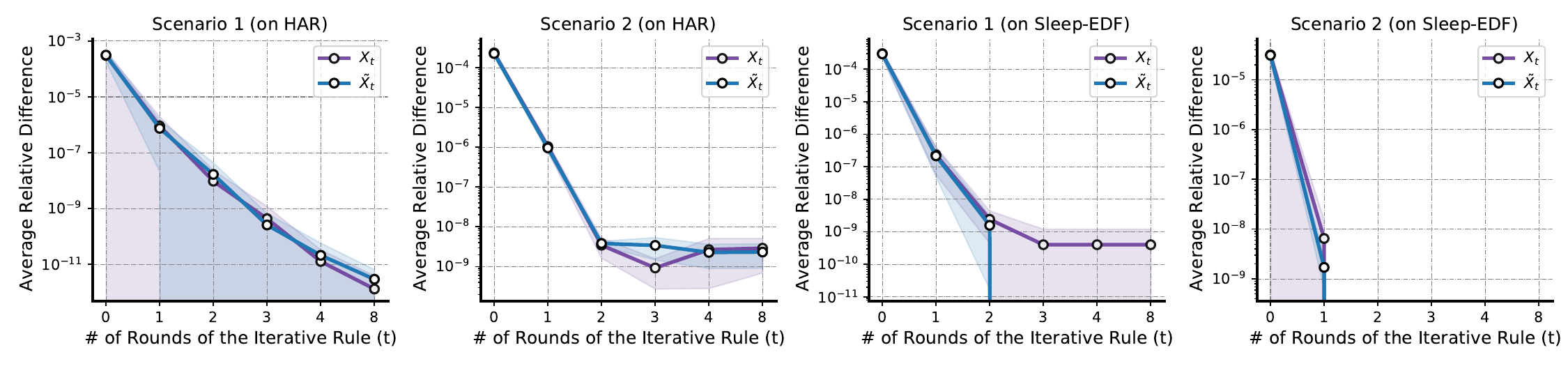}
    {\caption{Verification of Convergence Speed. When the \# of iteration $t$ becomes larger, the average relative difference can exceed the minimum precision and go to zero (like Scenario 2 on Sleep-EDF).}}
    \label{fig:linear_convergence}
\end{figure}

{\bf One Round is Enough.} 
We consider the performance of downstream classification with different rounds of iterative rules. {The results are shown in Table~3 and Table~4}.
\begin{table}[h!]\caption{Performance with Different Rounds (on HAR)}
\centering
\resizebox{0.9\textwidth}{!}
{
\begin{tabular}{c|ccccc} 
\toprule 
\# of rounds (t) & 1 & 2 & 3 & 4 & 8 \\
\midrule 
time per sweep & 8.774s & 10.316s & 11.224s & 12.781s & 17.402s \\
accuracy (\%) & 93.30 $\pm$ 0.413 & 93.35 $\pm$ 0.172 & 93.35 $\pm$ 0.150 & 93.35 $\pm$ 0.122 & 93.35 $\pm$ 0.122 \\

\bottomrule 
\end{tabular}}
\label{tb:number_of_recursion1} \end{table}

\begin{table}[h!]{\caption{Performance with Different Rounds (on Sleep-EDF)}}
\centering
\resizebox{0.9\textwidth}{!}
{
{\begin{tabular}{c|ccccc} 
\toprule 
\# of rounds (t) & 1 & 2 & 3 & 4 & 8 \\
\midrule 
time per sweep & 148.375s & 160.924s & 173.361s & 188.193s & 242.386s \\
accuracy (\%) & 85.25 $\pm$ 0.209 & 85.33 $\pm$ 0.173 & 85.31 $\pm$ 0.178 & 85.31 $\pm$ 0.177 & 85.31 $\pm$ 0.177\\
\bottomrule 
\end{tabular}}}
\label{tb:number_of_recursion2} \end{table}

We observe that with an increasing number of rounds of iterative rule, the classification results will not improve further, however, the time consumption increases. Thus, we use only one round of the iterative rule in our experiments.

\section{Algorithm}
For small tensors, the algorithmic pipeline is presented in Algorithm~\ref{algo:sao1}. For large tensors, we use mini-batch alternating least squares optimization in Algorithm~\ref{algo:sao}.
\begin{algorithm}[t!]\footnotesize
	\SetAlgoLined
	\textbf{Input:} Data tensor $\mathcal{T}\in\mathbb{R}^{N\times I\times J\times K}$; initialized $\mathbf{A},\mathbf{B},\mathbf{C}, \tilde{\mathbf{X}}, \mathbf{X}$; other hyperparameters $\alpha,\beta,\gamma$\;
	Obtain the augmented tensor $\tilde{\mathcal{T}}$\;
	\Repeat{max sweep exceeds or change of loss $<0.1\%$ within 3 consecutive sweeps}{
		Use $\mathbf{A},\mathbf{B},\mathbf{C},\tilde{\mathbf{X}}$ to update $\mathbf{X}$ by our iterative rules (one iteration) in Eqn.~\eqref{eq:improved_x}\;
		Use $\mathbf{A},\mathbf{B},\mathbf{C},\mathbf{X}$ to update $\tilde{\mathbf{X}}$ by our iterative rules (one iteration) in Eqn.~\eqref{eq:improved_x}\;
		Use $\mathbf{B},\mathbf{C},\mathbf{X}, \tilde{\mathbf{X}}$ to update $\mathbf{A}$ by solving least square problem\;
		Use $\mathbf{A},\mathbf{C},\mathbf{X}, \tilde{\mathbf{X}}$ to update $\mathbf{B}$ by solving least square problem\;
		Use $\mathbf{A},\mathbf{B},\mathbf{X}, \tilde{\mathbf{X}}$ to update $\mathbf{C}$ by solving least square problem\;
		}
	\textbf{Output:} the learned bases $\{\mathbf{A},\mathbf{B},\mathbf{C}\}$.
	\caption{Alternating Least Squares}
	\label{algo:sao1}
\end{algorithm}

\begin{algorithm}[t!]\footnotesize
	\SetAlgoLined
	\textbf{Input:} Data tensor $\mathcal{T}\in\mathbb{R}^{N\times I\times J\times K}$; initialized $\{\mathbf{A}^1,\mathbf{B}^1,\mathbf{C}^1\}$; batch size $b$; learning rate $\eta$; other hyperparameters $\alpha,\beta,\gamma$; initial counter $l=1$\;
	\Repeat{max sweep exceeds or change of loss $<0.1\%$ within 3 consecutive sweeps}{
		shuffle the data tensor $\mathcal{T}$; \quad $/*$ start a new sweep $*/$\\
		\For{a tensor batch $\mathcal{T}^l\in\mathbb{R}^{b\times I\times J\times K}$ and its augmentation $\tilde{\mathcal{T}}^l$}{
			{
			    use $\mathbf{A}^l,\mathbf{B}^l,\mathbf{C}^l$ to initialize $\mathbf{X}$ based on $\mathcal{T}^l$\;
			    use $\mathbf{A}^l,\mathbf{B}^l,\mathbf{C}^l$ to initialize $\tilde{\mathbf{X}}$ based on $\tilde{\mathcal{T}}^l$\;
				Use $\mathbf{A}^l,\mathbf{B}^l,\mathbf{C}^l,\tilde{\mathbf{X}}$ to update $\mathbf{X}$ by our iterative rules (one iteration) in Eqn.~\eqref{eq:improved_x}\;
		Use $\mathbf{A}^l,\mathbf{B}^l,\mathbf{C}^l,\mathbf{X}$ to update $\tilde{\mathbf{X}}$ by our iterative rules (one iteration) in Eqn.~\eqref{eq:improved_x}\;
		Use $\mathbf{B}^l,\mathbf{C}^l,\mathbf{X}, \tilde{\mathbf{X}}$ to obtain $\mathbf{A}^{l+1}$ by solving least square problem\;
		Use $\mathbf{A}^{l+1},\mathbf{C}^l,\mathbf{X}, \tilde{\mathbf{X}}$ to obtain $\mathbf{B}^{l+1}$ by solving least square problem\;
		Use $\mathbf{A}^{l+1},\mathbf{B}^{l+1},\mathbf{X}, \tilde{\mathbf{X}}$ to obtain $\mathbf{C}^{l+1}$ by solving least square problem\;
			}
			$l=l+1$ \quad $/*$ increment the counter $*/$\;}}
	\textbf{Output:} the learned bases $\{\mathbf{A}^L,\mathbf{B}^L,\mathbf{C}^L\}$.
	\caption{Mini-batch Alternating Least Squares}
	\label{algo:sao}
\end{algorithm}

\section{Experimental Details}

\subsection{Dataset Processing} 
{\em Sleep-EDF} is public, which contains 153 full-night EEG (from Fpz-Cz and Pz-Oz electrode locations), EOG (horizontal), and submental chin EMG recordings, under
Open Data Commons Attribution License v1.0 and {\em MGH Sleep} is provided by \citep{biswal2018expert}, where F3-M2, F4-M1, C3-M2, C4-M1, O1-M2, O2-M1 channels are used, containing 6,478 recordings. These two EEG datasets are processed in a similar way. First, the raw data are (long) recordings of each subject. On subject-level, these recordings are categorized into unlabeled and labeled sets by $90\%:10\%$. Then the labeled sets are further separated into training and test by $5\%:5\%$. Next, within each set (unlabeled, training, test), recordings are further segmented into disjoint 30-second-long periods, which are the data samples in our study. Each data sample is represented as a matrix, {\em channel by timestamp}, and they are associated with one of five sleep stages, Awake (W), Non-REM stage 1 (N1), Non-REM stage 2 (N2), Non-REM stage 3 (N3), and REM stage (R). {\em HAR} is also public, collected as 3-axial linear acceleration and 3-axial angular velocity at a constant rate of 50Hz by the embedded accelerometer and gyroscope. It has been randomly partitioned into 70\% and 30\%. We use 70\% as the unlabeled data (labels removed) and split the other part: 15\% as training data and 15\% as the test data. The license of this
dataset is included in their citation. {\em PTB-XL} is a public ECG dataset, which contains 21,837 clinical 12-lead ECGs (male: 11,379 and female: 10,458) of 10-second length with a sampling frequency of 500 Hz. We randomly split the dataset into unlabeled (remove the labels) and labeled sets by $90\%:10\%$; then the labeled sets are further separated into training and test by $5\%:5\%$. This dataset is under Open BSD 3.0.

All datasets are de-identified (e.g., no names, no locations), and there is no offensive content. All the labels are also provided along with the datasets.
The label distributions are shown below.

\begin{table}[h!] \label{tb:datastatistics2} 
\centering
{\caption{Class Label Distribution}
\resizebox{0.9\textwidth}{!}{\begin{tabular}{c|c} 
\toprule 
  {\bf Name} & Label Distribution \\
  \midrule
  {Sleep-EDF} & W: 68.8\%, N1: 5.2\%, N2: 16.6\%, N3: 3.2\%, R: 6.2\% \\
  {HAR} & Walk: 16.72\%, Walk upstairs: 14.99\%, Walk downstairs: 13.65\%, \\
  &Sit: 17.25\%, Stand: 18.51\%, Lay: 18.88\% \\
  {PTB-XL} & Male: 52.11\%, Female: 47.89\% \\
  {MGH Sleep} & W: 44.3\%, N1: 9.9\%, N2: 14.4\%, N3: 17.6\%, R: 13.8\% \\
\bottomrule
\end{tabular}}}
\end{table}

\subsection{STFT Transform} We find that directly using the raw data (spatial information) does not provide good results, even for the deep learning models. Thus, we take Short-Time Fourier Transforms (STFT) as a preprocessing step. From a single channel, we can extract both the amplitude and phase information, which is then stacked together as two different channels. After STFT, each data sample becomes a three-order tensor, {\em channel by frequency by timestamp}. The FFT size is 256 and hop length is 32 for Sleep-EDF; the FFT size is 64 and hop length is 2 for HAR; the FFT size is 256 and hop length is 64 for PTB-XL; and the FFT size is 512 and hop length is 128 for MGH.  We use these third-order tensors as final input data samples for all models.

\subsection{Data Augmentation} 
In our work, we build our feature extractor $\mathbf{f}(\cdot)$ from tensor decomposition tool, which may not be as expressive/flexible as deep neural networks in SSL, and thus we also ask the augmentation methods to allow a similar component-based representation for the perturbed data. Use EEG signals as an example, we consider jittering and bandpass filtering as two augmentation methods in the experiments, which perturb the signal frequency information and will not significantly change the low-rank structure of the data.

As mentioned in the main text, we consider three different augmentation methods: (i) {\em Jittering} adds additional perturbations to each sample. We consider both high and low-frequency noise on each channel independently. For high-frequency noise, we first generate a noisy sequence $\mathbf{s}$, which has the same length as the signal channel, and each element of $\mathbf{s}$ is i.i.d. sampled from a uniform distribution $U[-1,1]$. We then control the amplitude of $\mathbf{s}$ by the noisy degree $d\in\mathbb{R}$. Finally, we add the scaled noisy sequence $d\cdot\mathbf{s}$ to the channel. In the experiment, $d=0.05$ for Sleep-EDF, $d=0.002$ for HAR, $d=0.001$ for PTB-XL and $d=0.01$ for MGH. For low-frequency noise, we generate a short noisy sequence (the length is randomly sampled from a uniform distribution $U[100, \mbox{length of channel}]$) in the same way and then use {\em scipy.interpolate.interp1d} to interpolate the noisy sequence to be at the same length as the channel. The choice of high-frequency noise or low-frequency noise, or both are coin-tossed with equal probability. (ii) {\em Bandpass filtering} reduces signal noise. We use the order-1 Butterworth filter by {\em scipy.signal.butter} to preserve only the within-band frequency information. The high-pass and low-pass are $(1Hz,30Hz)$ and $(10Hz,49Hz)$ for Sleep-EDF, $(1Hz,20Hz)$ and $(5Hz,24.5Hz)$ for HAR, $(1Hz,30Hz)$ and $(10Hz,50Hz)$ for PTB-XL, $(1Hz,30Hz)$ and $(10Hz,50Hz)$ for MGH. Low-pass or high-pass or both are selected with equal probability. Also, the bandpass filtering is applied to each channel independently. The intuition is that the low-pass signals and high-pass signals might be both useful. (iii) {\em 3D position rotation} is an augmentation technique used only for HAR datasets, which have x-y-z axis information from accelerometer and gyroscope sensors. We apply a 3D x-y-z coordinate system rotation by a rotation matrix to mimic different cellphone positions. 
All augmentation methods are applied in sequence (i) (ii) (iii). The STFT is performed after the data augmentation.

\begin{figure}[h!] 
    \centering
    \includegraphics[width=0.9\textwidth]{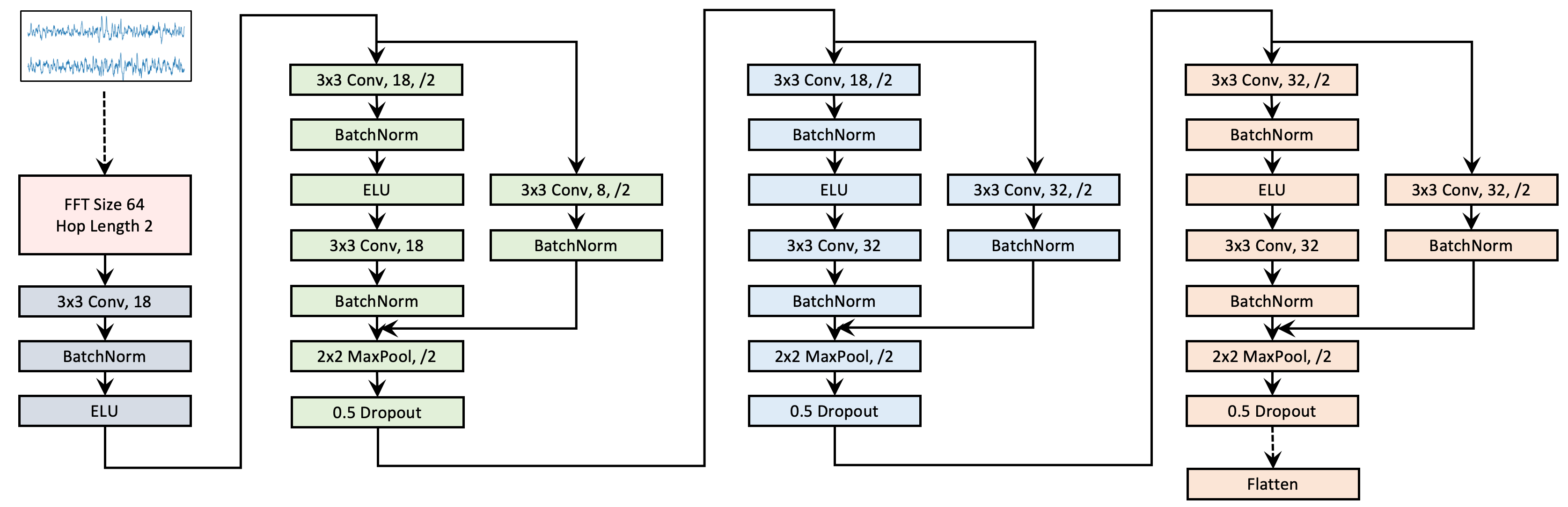}
    \caption{Backbone CNN Architecture (HAR). The architectures for other three datasets are similar.}
    \label{fig:backbone_cnn}
\end{figure}

\subsection{Implementation} Since all deep learning baselines are based on CNN, we use the same backbone model, shown in Figure~\ref{fig:backbone_cnn}. The model is adopted from \citep{cheng2020subject}. Based on the backbone model, we add a fully connected layer for the reference CNN model, add non-linear layers for self-supervised models (they also have their respective loss), and add corresponding deconvolutional layers for autoencoder models. The reference CNN models 
is end-to-end, and they are trained on the training set; other baselines learn a low-dimensional feature extractor from an unlabeled set, and then a logistic classifier is trained on the training set, on top of the feature extractor. Without loss of generality, we use $128$ as batch size. For deep learning models, we use Adam optimizer with a learning rate $1\times 10^{-3}$ and weight decaying $5\times 10^{-4}$. We use $2\times 10^{-3}$ as the learning rate for our \method. Since the paper deals with unsupervised learning and uses standard logistic regression ({\em sklearn.linear\_model.LogisticRegression}) to evaluate, the hyperparameters of all models (except supervised model) are chosen based on the classification performance on the training data. For our ATD model, by default, we use $R=32,\alpha=1\times 10^{-3},\beta=2$, $\gamma=b$ (which is the batch size), and all the reference configurations for each experiments are listed in code appendix. {The choice of $R$ depends on the trade-off between model fitness and time complexity. Specifically, a larger $R$ means better fitness and more preserved information in the extracted representations, while the number of learnable parameters and time complexity also increases linearly with $R$. In our paper, we run simple CP decomposition on a small subset of the tensor and monitor the fitness curve. We find that the fitness does not improve much around $R=32$ for all datasets (which means the real tensor might have a smaller rank and the residual part might be just noise). Thus, we choose $R=32$ throughout the paper.}

\subsection{Ablation Studies on  Data Augmentation}
\paragraph{Effect of High-quality Data Augmentations.} We study the effects of high-quality data augmentation methods. Let us assume {\em (a): jittering}, {\em (b): bandpass filtering}, {\em (c): 3D position rotation}. We test on different combinations of the augmentation methods. The experiments are conducted on two datasets: (i) for the MGH dataset, we use 50,000 unlabeled data, 5,000 training samples, and all test samples; (ii) for the HAR dataset, we use all data. We use 128 as the batch size for MGH Sleep, 64 for HAR, and $R=32$ as the tensor rank.
\begin{table}[h!]
\centering
\caption{Ablation Studies on Data Augmentation (MGH Sleep)}
\resizebox{0.9\textwidth}{!}{\begin{tabular}{cc|cc|cc} 
\toprule 
  {\bf Method} & {\bf Acc (\%)} & {\bf Method} & {\bf Acc (\%)} & {\bf Method} & {\bf Acc (\%)} \\
  \midrule
  (a) & 72.53 $\pm$ 0.652 & (b) & 73.60 $\pm$ 0.270 & (a)+(b) & 74.15 $\pm$ 0.203\\
\bottomrule
\end{tabular}} \label{tb:effect_data_augmentation_mgh} 
\end{table}

\begin{table}[h!]
\centering
\caption{Ablation Studies on Data Augmentation (HAR)}
\resizebox{0.9\textwidth}{!}{\begin{tabular}{cc|cc|cc} 
\toprule 
  {\bf Method} & {\bf Acc (\%)} & {\bf Method} & {\bf Acc (\%)} & {\bf Method} & {\bf Acc (\%)} \\
  \midrule
  (a) & 92.06 $\pm$ 0.621 & (a)+(b) & 92.70 $\pm$ 0.266 & (a)+(b)+(c) & 93.35 $\pm$ 0.357\\
  (b) & 91.99 $\pm$ 0.274 & (a)+(c) & 92.93 $\pm$ 0.590 &  / & /\\
  (c) & 92.29 $\pm$ 0.330 & (b)+(c) & 92.26 $\pm$ 0.479 &  / & /\\
\bottomrule
\end{tabular}} \label{tb:effect_data_augmentation_har} 
\end{table}

Table~\ref{tb:effect_data_augmentation_mgh} and Table~\ref{tb:effect_data_augmentation_har} conclude that the augmentation methods influence the final classification results. However, for different datasets, the effects are different. We observe that for the MGH Sleep dataset, bandpass filtering works better than jittering. In HAR, combinations of augmentations work better than the individual augmentations. 
% The reason might be that the changes brought by individual methods are still within the range of the class, while some super-impositions might change the class identity of the data sample, which hurts the classification. 
Overall, we find that with more diverse data augmentation methods, the final results are relatively better. The study of how to choose/design (or even automatically generate) better augmentation techniques will be our future work.

{\paragraph{Impact of Low-quality Data Augmentations.} 
We also study the impact of low-quality data augmentations. We use the SALS model (mini-batch tensor baseline without data augmentation) and our \method as the reference models and conduct the following experiments on MGH and HAR:
\begin{itemize}
    \item MGH-d: we change the $d$ values for the jittering data augmentation, $d$ means the ratio of the amplitude of the high/low frequency noise over the amplitude of the signal.
    \item MGH-(A): on MGH data, we set the degree of the {\em jittering} method to an unrealistic value $d=1$ as the low-quality augmentation method, meaning that the magnitude of the noise is the same as the magnitude of the signals. We keep the {\em bandpass filtering} unchanged.
    \item MGH-(B): on MGH data, we randomly create a {\em bandpass filter} from $(1Hz, 10Hz)$ or $(30Hz, 50Hz)$ as the low-quality augmentation method, which drops the critical middle-band information. We keep the {\em jittering} unchanged.
    \item MGH-(C): on MGH data, we combine the above two low-quality augmentation methods.
    \item HAR-(A)(B)(C): follow the similar low-quality data augmentation method design on HAR. For (B), we use $(1Hz, 5Hz)$ and $(20Hz, 24.5Hz)$ as the low-quality bandpass.
\end{itemize}

The performances are shown in Table~8. We find that low-quality data augmentations will hurt the learned tensor bases and hinder the downstream classification. With $d$ changing from $0.02$ to $5$, the {\em jittering} method becomes more unrealistic and the generated samples deviate from the real signal data distribution. Thus, we can find that the final classification performance also becomes worse gradually. If all data augmentation methods are of low-quality, the performance cannot surpass the base SALS model. We also find that the performance drop is not significant. The reason might be that even the augmentation methods are of low-quality, the Frobenius loss can still enforce the model to learn decent subspaces for better fitness and find generic low-rank features (opposed to the classification-oriented low-rank features), in this case, the features only follow fitness principle not the alignment principle, so the performance will be similar compared to SALS.
}
\begin{table*}[t!] \centering
	{\caption{Performance of Low-quality Data Augmentation (on MGH and HAR)}}
	\resizebox{0.95\textwidth}{!}{
		{\begin{tabular}{c|cccccc} \toprule 
		MGH & $d=0.02$ & $d=0.05$ & $d=0.1$ & $d=0.5$ & $d=5$ \\
	    Accuracy & 74.18 $\pm$ 0.326 & 74.10 $\pm$ 0.302 & 73.85 $\pm$ 0.530 & 73.39 $\pm$ 0.493 & 72.18 $\pm$ 0.676 \\
		\midrule
		MGH & SALS & \method & (A) & (B) & (C) \\
	    Accuracy & 71.93 $\pm$ 0.379 & 74.15 $\pm$ 0.431 & 72.73 $\pm$ 0.624 & 72.10 $\pm$ 0.719 & 70.75 $\pm$ 0.771 \\
	    \midrule
		HAR & SALS & \method & (A) & (B) & (C) \\
	    Accuracy & 91.86 $\pm$ 0.295 & 93.35 $\pm$ 0.357 & 92.04 $\pm$ 0.308 & 92.48 $\pm$ 0.469 & 91.43 $\pm$ 0.835 \\
		\bottomrule
	\end{tabular}}\label{tb:low-quality-data-augmentation}}
\end{table*}

\subsection{Ablation Studies on Hyperparameters}
This section conducts ablation studies for decomposition rank $R$ and other hyperparameters, $\alpha$, $\beta$, $\gamma$. The experiments are conducted on Sleep-EDF with 50,000 random unlabeled data, 5,000 random training samples, and all test samples, and the HAR dataset.

\begin{figure}[t!] 
    \centering
    \includegraphics[width=0.9\textwidth]{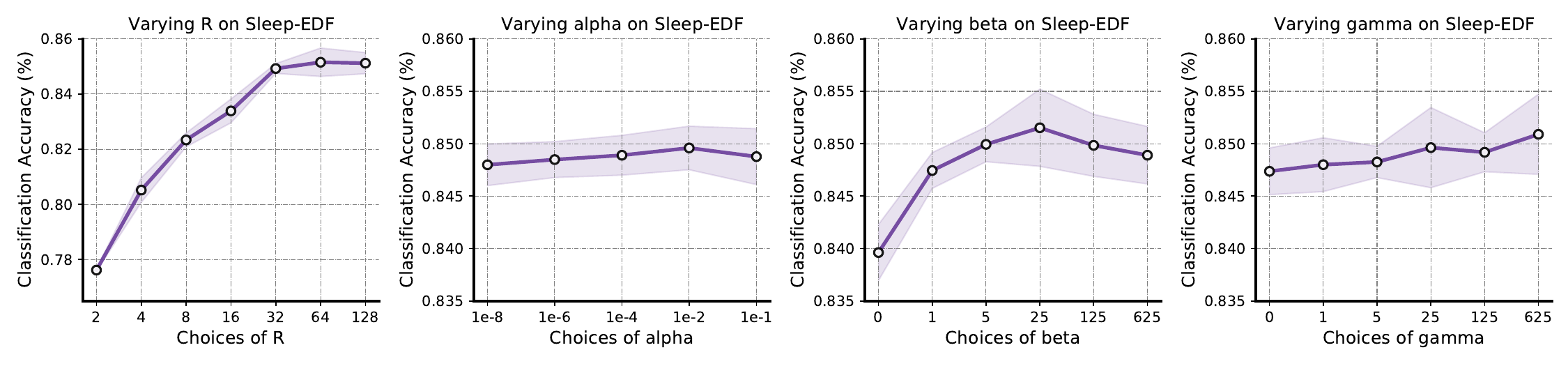}
    \includegraphics[width=0.9\textwidth]{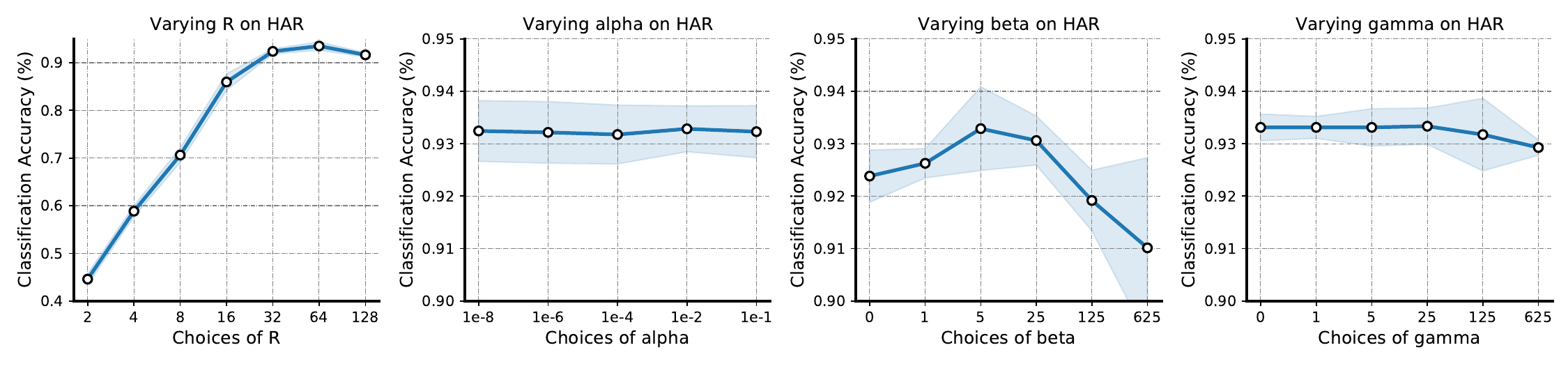}
    \caption{Ablation Studies on Hyperparameters}
    \label{fig:hyperparameter_sleep}
\end{figure}

\begin{table*}[t!] \centering
	{\caption{Comparison with Tensor-based Methods with Different R (on HAR)}}
	\resizebox{0.95\textwidth}{!}{
		{\begin{tabular}{c|ccccc} \toprule 
		{\bf Model} & {\bf $R=8$} & {\bf $R=16$} & {\bf $R=32$} & {\bf $R=64$} & {\bf $R=128$} \\
		\midrule
        SALS & 69.59 $\pm$ 0.526 & 83.92 $\pm$ 0.416 & 91.84 $\pm$ 0.295 & 91.89 $\pm$ 0.217 & 91.55 $\pm$ 0.388 \\
        GR-SALS & 69.62 $\pm$ 0.458 & 84.20 $\pm$ 0.727 & 92.33 $\pm$ 0.282 & 92.28 $\pm$ 0.359 & 91.84 $\pm$ 0.534 \\
        $\method_{ss-}$ & 70.27 $\pm$ 0.488 & 84.84 $\pm$ 0.557 & 92.41 $\pm$ 0.391 & 92.71 $\pm$ 0.243 & 92.32 $\pm$ 0.330 \\
        \method & 71.91 $\pm$ 0.253 & 85.61 $\pm$ 0.294 & 93.35 $\pm$ 0.357 & 93.43 $\pm$ 0.411 & 92.97 $\pm$ 0.273 \\
		\bottomrule
	\end{tabular}}\label{tb:varying_R_comparison_HAR}}
\end{table*}
\begin{table*}[t!] \centering
	{\caption{Comparison with Tensor-based Methods with Different R (on Sleep-EDF)}}
	\resizebox{0.95\textwidth}{!}{
		{\begin{tabular}{c|ccccc} \toprule 
		{\bf Model} & {\bf $R=8$} & {\bf $R=16$} & {\bf $R=32$} & {\bf $R=64$} & {\bf $R=128$} \\
		\midrule
        SALS & 81.26 $\pm$ 0.345 & 82.59 $\pm$ 0.638 & 84.27 $\pm$ 0.481 & 84.55 $\pm$ 0.527 & 84.49 $\pm$ 0.317 \\
        GR-SALS & 81.72 $\pm$ 0.664 & 82.74 $\pm$ 0.481 & 84.33 $\pm$ 0.356 & 84.87 $\pm$ 0.486 & 84.90 $\pm$ 0.781 \\
        $\method_{ss-}$ & 81.27 $\pm$ 0.568 & 82.5 $\pm$ 0.674 & 84.19 $\pm$ 0.221 & 84.47 $\pm$ 0.258 & 84.44 $\pm$ 0.577 \\
        \method & 82.49 $\pm$ 0.464 & 83.31 $\pm$ 0.591 & 85.01 $\pm$ 0.224 & 85.30 $\pm$ 0.483 & 85.32 $\pm$ 0.305 \\
		\bottomrule
	\end{tabular}}\label{tb:varying_R_comparison_sleep}}
\end{table*}

The results are shown in Figure~\ref{fig:hyperparameter_sleep}. First, we can conclude that with a larger decomposition rank $R$, the performance will be better generally. Though we observe that the performance worsens from $R=64$ to $R=128$, with limited training data, {if the representation size (equals to $R$) becomes larger, the logistic regression model can overfit. 
We compare our model to other tensor-based methods with $R=8, 16, 32, 64, 128$ and the Table~9 and Table~10 shows that \method outperforms the tensor baselines consistently with different R.}

Also, we find that the choices of $\alpha$ do not affect the final performance a lot. {Finally, we find that it is
easy for users to select the $\beta$ values from a large range in general. For example, 
selecting a $\beta\in[5, 125]$ would guarantee good results on Sleep-EDF while on
HAR, the selection range is $[5, 25]$. The choice of $\beta$ does affect the
final performance, but it is not tricky to search for a $\beta$ for good performance. We also find that the accuracy score first increases then decreases with an increasing value of $\beta$. The reason might be that a large $\beta$ will negatively affect the fitness loss.}

{\subsection{Statistical Testing and Running Time Comparison}
In this section, we conduct T-test on the result in Section~4.2.1 and calculate the p-values
in the parenthesis of Table~\ref{tb:time_table1} and Table~\ref{tb:time_table2} (the experimental results are copied from Table~2). Commonly, a p-value smaller than 0.05
would be considered as significant. We can see that our model show significant performance gain over all baselines on all tasks. 

We have also reported the running time per sweep/epoch in the tables. When recording the running time, we duplicated the environment mentioned in Section~4.1, stopped other programs and ran all the models one by one on GPUs. We record the first 8 sweeps/epochs of all models and drop the first 3 sweeps (since they might be unstable). The average running time of the last 5 epochs are reported in Table~\ref{tb:time_table1} and Table~\ref{tb:time_table2} while the accuracy results are from Table~2. Note that on HAR, PTB-XL and Sleep-EDF, all methods use 128 as the batch size and on MGH, all methods use 512 as the batch size. The tensor based methods all use $R=32$ as the rank. We can conclude that the tensor based methods are generally more time-efficient than the deep learning methods with fewer parameters. Since our model \method and the variant $\method_{ss-}$ use the augmented tensors, they cost more compared to other tensor based methods (since the size of training tensors doubles), however, we also observe empirically that they can converge faster with around half number of the epochs.
}

\begin{table*}[t!] \centering
	\caption{Result Significance and Running Time (\%) for Sleep-EDF and HAR}
	\resizebox{0.9\textwidth}{!}{
		\begin{tabular}{lcccccc} \toprule 
			\multirow{2}{*}{} & \multicolumn{3}{c}{\bf Sleep-EDF (5,000)} & \multicolumn{3}{c}{\bf HAR (1,473)}\\
			\cmidrule(r){2-4}  \cmidrule(r){5-7}
			& {Accuracy} & {\# of Params.} & {Time per sweep} & {Accuracy} & {\# of Params.} & {Time per sweep}\\
			\cmidrule{1-7} 
			{\bf Self-sup models:} \\
			SimCLR-32 & 84.98 $\pm$ 0.358 & 210,384 & {260.299s} & 74.75 $\pm$ 0.723 & 53,286 & {8.459s} \\
			SimCLR-128 & {\bf 85.19 $\pm$ 0.358} & 222,768 & {265.809s} & 76.69 $\pm$ 0.697 & 65,670 & {8.532s}\\
			BYOL-32 & 84.29 $\pm$ 0.405 & 211,440 & {255.614s} & 73.71 $\pm$ 2.832 & 54,342 & {8.430s}\\
			BYOL-128 & 83.26 $\pm$ 0.337 & 239,280 & {257.266s} & 71.79 $\pm$ 1.866 & 82,182 & {8.478s}\\
			\midrule
			{\bf Auto-encoders:} \\
			AE-32 & 74.78 $\pm$ 0.723 & 217,216 & {153.684s} & 63.13 $\pm$ 0.775 & 62,940 & {7.530s}\\
			AE-128 & 75.17 $\pm$ 0.897 & 241,888 & {156.813s} & 60.52 $\pm$ 1.604 & 87,612 & {7.662s}\\
			$\mbox{AE}_{ss}$-32 & 80.92 $\pm$ 0.345 & 217,216 & {301.773s} & 71.70 $\pm$ 2.135 & 62,940 & {7.765s}\\
			$\mbox{AE}_{ss}$-128 & 81.84 $\pm$ 0.259 & 241,888 & {307.546s} & 72.43 $\pm$ 1.370 & 87,612 & {7.804s}\\
			\midrule
			{\bf Tensor models:} \\
			SALS & 84.27 $\pm$ 0.481 {(0.0041)} & 7,328 & {86.281s} & 91.86 $\pm$ 0.295 {(2e-5)} & 2,688 & {7.535s}\\
			GR-SALS & 84.33 $\pm$ 0.356 {(0.0019)} & 7,328 & {109.916s} &  92.33 $\pm$ 0.282 {(0.0003)} & 2,688 & {7.829s}\\
			$\method_{ss-}$ & 84.19 $\pm$ 0.221 {(9e-5)} & 7,328 & {147.568s} & 92.41 $\pm$ 0.391 {(0.0011)} & 2,688 & {8.604s}\\
			\method & {85.01 $\pm$ 0.224} & 7,328 & {148.375s} & {\bf 93.35 $\pm$ 0.357} & 2,688 & {8.672s}\\
			\bottomrule
			\multicolumn{5}{l}{result format: mean $\pm$ standard deviation (p-value)}
	\end{tabular}\label{tb:time_table1}}
\end{table*}

\begin{table*}[t!] \centering
	\caption{Result Significance and Running Time (\%) for PTB-XL and MGH}
	\resizebox{0.9\textwidth}{!}{
		\begin{tabular}{lcccccc} \toprule 
			\multirow{2}{*}{} & \multicolumn{3}{c}{\bf PTB-XL (2,183)} & \multicolumn{3}{c}{\bf MGH (5,000)}\\
			\cmidrule(r){2-4}  \cmidrule(r){5-7}
			& {Accuracy} & {\# of Params.} & {Time per sweep} & {Accuracy} & {\# of Params.} & {Time per sweep}\\
			\cmidrule{1-7} 
			{\bf Self-sup models:} \\
			%   MoCo & 64.65 $\pm$ 0.615 & 74.99 $\pm$ 0.713 & 74,626 & 65.78 $\pm$ 1.914 & 70.51 $\pm$ 0.647 & 172,104\\
			SimCLR-32 & 69.25 $\pm$ 0.355 & 200,960 & {18.714s} & 67.34 $\pm$ 0.970 & 212,624 & {1449.368s}\\
			SimCLR-128 & 68.19 $\pm$ 0.793 & 237,920 & {19.037s} & 66.98 $\pm$ 1.331 & 246,608 & {1457.283s}\\
			BYOL-32 & 65.08 $\pm$ 1.535 & 202,016  & {18.410s} & 68.83 $\pm$ 1.168 & 214,736 & {1451.468s}\\
			BYOL-128 & 65.49 $\pm$ 0.612 & 254,432 & {18.680s} & 68.55 $\pm$ 1.339 & 279,632 & {1461.181s}\\
			\midrule
			{\bf Auto-encoders:} \\
			AE-32 & 59.01 $\pm$ 0.896  & 224,528 & {11.229s} & 68.58 $\pm$ 0.427 & 220,088 & {851.118s}\\
			AE-128 & 58.29 $\pm$ 0.412 & 298,352 & {11.396s} & 67.05 $\pm$ 1.375 & 257,048 & {815.858s}\\
			$\mbox{AE}_{ss}$-32  & 68.47 $\pm$ 0.231 & 224,528 & {18.263s} & 71.46 $\pm$ 0.386 &  220,088 & {1486.244s}\\
			$\mbox{AE}_{ss}$-128 & 68.88 $\pm$ 0.604 & 298,352 & {18.465s} & 70.19 $\pm$ 0.617 & 257,048 & {1504.545s}\\
			\midrule
			{\bf Tensor models:} \\
			SALS & 69.15 $\pm$ 0.483 {(0.0023)} & 7,296 & {8.988s} & 71.93 $\pm$ 0.379 {(5e-6)} & 9,984 & {782.763s}\\
			GR-SALS & 69.02 $\pm$ 0.477 {(0.0012)} & 7,296 & {9.747s} & 72.35 $\pm$ 0.228 {(8e-6)} & 9,984 & {970.292s}\\
			$\method_{ss-}$ & 69.38 $\pm$ 0.612 {(0.0129)} & 7,296 & {12.560s} & 72.78 $\pm$ 0.522 {(0.0005)} & 9,984 & {1327.188s}\\
			\method & {\bf 70.26 $\pm$ 0.523} & 7,296 & {12.599s} & {\bf 74.15 $\pm$ 0.431} & 9,984 & {1360.569s}\\
			\bottomrule
			\multicolumn{5}{l}{result format: mean $\pm$ standard deviation (p-value)}
	\end{tabular}\label{tb:time_table2}}
\end{table*}

{\subsection{Comparison with Supervised Tensor Learning}
In this subsection, we compare our model with two supervised tensor learning baselines, UMLDA \cite{lu2008uncorrelated} and supervised tensor learning (STL) \cite{tao2005supervised}. {\bf UMLDA} extracts uncorrelated discriminative features by sequential tensor-to-vector projections, and it includes two stages: in the first stage (supervised pre-training stage), it uses label information to maximize the Fisher’s discrimination criterion (FDC) as the objective and extract uncorrelated features (i.e., representations); in the second stage (supervised learning stage), it uses another supervised model to map the representations to the labels. To make a fair comparison, we use $32$ as the uncorrelated feature dimension (thus, it has the same number of learnable parameters as our model) and also use logistic regression (LR) for the second stage. {\bf STL} is an end-to-end supervised tensor learning baseline, and it is originally proposed for binary classification with a rank-one parameterized tensor (i.e., outer product of multiple vectors). In the comparison, we extend STL for mult-class classification by including more rank-one parameterized tensors (one for each class). We use cross entropy loss to optimize the revised STL model. Both baselines are implemented with PyTorch and use the training and test set only. We have carefully turned the baseline models to achieve higher accuracy. 

\paragraph{Result Analysis.} The results are shown in Table~13. We can conclude that our methods outperform these two supervised tensor learning baselines significantly. The performance gap between UMLDA and our model can be explained by (i) UMLDA uses FDC criterion to design the loss function and also forces the learned feature dimensions to be uncorrelated, which might discard some essential class-dependent information; (ii) our model utilizes a large set of unlabeled data. STL gives poor performance because it is essentially a multilinear method with much fewer parameters than UMLDA and our \method, which hurts the expressive.

\begin{table*}[t!] \centering
	{\caption{Comparison with Supervised Tensor Learning}}
	\resizebox{1.0\textwidth}{!}{
		{\begin{tabular}{c|cccc} \toprule 
		{\bf Model} & {\bf Sleep-EDF (5,000)} & {\bf HAR (1,473)} & {\bf PTB-XL (2,183)} & {\bf MGH (5,000)} \\
		\midrule
		UMLDA (supervised pretrain + supervised LR) & 81.06 $\pm$ 0.093 & 85.73 $\pm$ 1.169 & 65.55 $\pm$ 0.267 & 62.04 $\pm$ 0.722\\
		STL (end-to-end supervised) & 77.86 $\pm$ 0.816 & 80.52 $\pm$ 0.189 & 61.83 $\pm$ 0.712 & 41.44 $\pm$ 0.597\\
		\method (unsupervised pretrain + supervised LR) & {\bf 85.01 $\pm$ 0.224} & {\bf 93.35 $\pm$ 0.357} & {\bf 70.26 $\pm$ 0.523} & {\bf 74.15 $\pm$ 0.431} \\
		\bottomrule
	\end{tabular}}\label{tb:supervised_comparison}}
\end{table*}
}
{\section{Derivation of Two-sided Bound}
We recall the definition of $\mathcal{L}_{ss}$ and $\mathcal{L}^{\Theta}_{ss}$ from Section~3.1.
\begin{align}
	\mathcal{L}_{ss} =&~~\mathcal{L}_{pos} + \lambda\mathcal{L}_{neg} \notag\\
	=&~ \mathbb{E} \left[\frac{\lambda}{1-r_p}\mbox{sim}\left(\mathbf{f}\left(\mathcal{X}_p \right),\mathbf{f}\left({\mathcal{Y}_q}\right)\right)\right] - \mathbb{E}\left[\left(\frac{\lambda r_p}{1-r_p}+1\right)\mbox{sim}\left(\mathbf{f}\left(\mathcal{X}_p \right),\mathbf{f}\left({\mathcal{Y}_q}\right)\right) \mid p = q \right], \notag\\
	\mathcal{L}^{\Theta}_{ss}(\gamma) =&~(\gamma+1)\mathbb{E} \left[\mbox{sim}\left(\mathbf{f}\left(\mathcal{X}_p\right),\mathbf{f}\left({\mathcal{Y}_q}\right)\right)\right] -\mathbb{E} \left[\mbox{sim}\left(\mathbf{f}\left(\mathcal{X}_p\right),\mathbf{f}\left({\mathcal{Y}_q}\right)\right)\mid p= q\right], \notag
\end{align}
and we want to prove that 
\begin{equation}
	C_1\mathcal{L}^{\Theta}_{ss}\left(\frac{\lambda-1}{C_1}\right) \leq \mathcal{L}_{ss} \leq C_2\mathcal{L}^{\Theta}_{ss}\left(\frac{\lambda-1}{C_2}\right),~
C_1=1+\max_{p}\frac{\lambda r_p}{1-r_p},C_2=1+\min_{p}\frac{\lambda r_p}{1-r_p}, \notag
\end{equation}
where $r_p$ is the label rate of class-$p$.

\begin{proof}
We start by arranging $\mathcal{L}_{ss}$,
\begin{align}
    \mathcal{L}_{ss}
	=&~ \mathbb{E} \left[\frac{\lambda}{1-r_p}\mbox{sim}\left(\mathbf{f}\left(\mathcal{X}_p \right),\mathbf{f}\left({\mathcal{Y}_q}\right)\right)\right] - \mathbb{E}\left[\left(\frac{\lambda r_p}{1-r_p}+1\right)\mbox{sim}\left(\mathbf{f}\left(\mathcal{X}_p \right),\mathbf{f}\left({\mathcal{Y}_q}\right)\right) \mid p = q \right] \notag\\
	=&~ \mathbb{E} \left[\left(\frac{\lambda r_p}{1-r_p} + \lambda\right)\mbox{sim}\left(\mathbf{f}\left(\mathcal{X}_p \right),\mathbf{f}\left({\mathcal{Y}_q}\right)\right)\right] - \mathbb{E}\left[\left(\frac{\lambda r_p}{1-r_p}+1\right)\mbox{sim}\left(\mathbf{f}\left(\mathcal{X}_p \right),\mathbf{f}\left({\mathcal{Y}_q}\right)\right) \mid p = q \right]\notag \\
	=&~ \mathbb{E} \left[\left(\frac{\lambda r_p}{1-r_p} + 1\right)\mbox{sim}\left(\mathbf{f}\left(\mathcal{X}_p \right),\mathbf{f}\left({\mathcal{Y}_q}\right)\right)\right]
	+ \mathbb{E} \left[\left(\lambda- 1\right)\mbox{sim}\left(\mathbf{f}\left(\mathcal{X}_p \right),\mathbf{f}\left({\mathcal{Y}_q}\right)\right)\right] \notag\\
	&- \mathbb{E}\left[\left(\frac{\lambda r_p}{1-r_p}+1\right)\mbox{sim}\left(\mathbf{f}\left(\mathcal{X}_p \right),\mathbf{f}\left({\mathcal{Y}_q}\right)\right) \mid p = q \right] \label{eq:exp_split1}\\
    =&~ \mathbb{E}_p\mathbb{E}_{q,\mathcal{X}_p,\mathcal{Y}_q} \left[\left(\frac{\lambda r_p}{1-r_p} + 1\right)\mbox{sim}\left(\mathbf{f}\left(\mathcal{X}_p \right),\mathbf{f}\left({\mathcal{Y}_q}\right)\right)\right]
	+ \mathbb{E} \left[\left(\lambda- 1\right)\mbox{sim}\left(\mathbf{f}\left(\mathcal{X}_p \right),\mathbf{f}\left({\mathcal{Y}_q}\right)\right)\right] \notag\\
	&- \mathbb{E}_p\mathbb{E}_{q,\mathcal{X}_p,\mathcal{Y}_q}\left[\left(\frac{\lambda r_p}{1-r_p}+1\right)\mbox{sim}\left(\mathbf{f}\left(\mathcal{X}_p \right),\mathbf{f}\left({\mathcal{Y}_q}\right)\right) \mid p = q \right] \label{eq:exp_split2}\\
	=&~ \mathbb{E}_p\left[\left(\frac{\lambda r_p}{1-r_p} + 1\right)\mathbb{E}_{q,\mathcal{X}_p,\mathcal{Y}_q} \left[\mbox{sim}\left(\mathbf{f}\left(\mathcal{X}_p \right),\mathbf{f}\left({\mathcal{Y}_q}\right)\right)\right]\right]
	+ \mathbb{E} \left[\left(\lambda- 1\right)\mbox{sim}\left(\mathbf{f}\left(\mathcal{X}_p \right),\mathbf{f}\left({\mathcal{Y}_q}\right)\right)\right] \notag\\
	&- \mathbb{E}_p\left[\left(\frac{\lambda r_p}{1-r_p}+1\right)\mathbb{E}_{q,\mathcal{X}_p,\mathcal{Y}_q}\left[\mbox{sim}\left(\mathbf{f}\left(\mathcal{X}_p \right),\mathbf{f}\left({\mathcal{Y}_q}\right)\right) \mid p = q \right]\right] \notag\\
	=&~ \mathbb{E}_p\left[\left(\frac{\lambda r_p}{1-r_p} + 1\right)\left(\mathbb{E}_{q,\mathcal{X}_p,\mathcal{Y}_q} \left[\mbox{sim}\left(\mathbf{f}\left(\mathcal{X}_p \right),\mathbf{f}\left({\mathcal{Y}_q}\right)\right)\right]-\mathbb{E}_{q,\mathcal{X}_p,\mathcal{Y}_q}\left[\mbox{sim}\left(\mathbf{f}\left(\mathcal{X}_p \right),\mathbf{f}\left({\mathcal{Y}_q}\right)\right) \mid p = q \right]\right)\right] \notag\\
	& + \mathbb{E} \left[\left(\lambda- 1\right)\mbox{sim}\left(\mathbf{f}\left(\mathcal{X}_p \right),\mathbf{f}\left({\mathcal{Y}_q}\right)\right)\right] \label{eq:inequality_1}\\
	\leq&~ \mathbb{E}_p\left[C_2\left(\mathbb{E}_{q,\mathcal{X}_p,\mathcal{Y}_q} \left[\mbox{sim}\left(\mathbf{f}\left(\mathcal{X}_p \right),\mathbf{f}\left({\mathcal{Y}_q}\right)\right)\right]-\mathbb{E}_{q,\mathcal{X}_p,\mathcal{Y}_q}\left[\mbox{sim}\left(\mathbf{f}\left(\mathcal{X}_p \right),\mathbf{f}\left({\mathcal{Y}_q}\right)\right) \mid p = q \right]\right)\right] \notag\\
	& + \mathbb{E} \left[\left(\lambda- 1\right)\mbox{sim}\left(\mathbf{f}\left(\mathcal{X}_p \right),\mathbf{f}\left({\mathcal{Y}_q}\right)\right)\right] \label{eq:inequality_2} \\
	=&~ (C_2+\lambda-1)\mathbb{E} \left[\mbox{sim}\left(\mathbf{f}\left(\mathcal{X}_p \right),\mathbf{f}\left({\mathcal{Y}_q}\right)\right)\right]-C_2\mathbb{E}\left[\mbox{sim}\left(\mathbf{f}\left(\mathcal{X}_p \right),\mathbf{f}\left({\mathcal{Y}_q}\right)\right) \mid p = q \right] \notag\\
	=&~ C_2\mathcal{L}^{\Theta}_{ss}\left(\frac{\lambda-1}{C_2}\right). \notag
\end{align}

From Eqn.~\eqref{eq:exp_split1} to Eqn.~\eqref{eq:exp_split2}, we use the fact that ``$\mathbb{E}[\cdot]$ means the expectation is taken over four interdependent random variables, i.e., $p,q,\mathcal{X}_p,\mathcal{Y}_q$", which is mentioned in Section~3.1. From  Eqn.~\eqref{eq:inequality_1} to Eqn.~\eqref{eq:inequality_2}, we use the fact that given $p$, the similarity of random pairs is smaller than the similarity of positive pairs $\mathbb{E}_{q,\mathcal{X}_p,\mathcal{Y}_q} \left[\mbox{sim}\left(\mathbf{f}\left(\mathcal{X}_p \right),\mathbf{f}\left({\mathcal{Y}_q}\right)\right)\right]\leq\mathbb{E}_{q,\mathcal{X}_p,\mathcal{Y}_q}\left[\mbox{sim}\left(\mathbf{f}\left(\mathcal{X}_p \right),\mathbf{f}\left({\mathcal{Y}_q}\right)\right) \mid p = q \right]$. The upper bound is derived by replacing $\frac{\lambda r_p}{1-r_p} + 1,~\forall p$ with $C_2=1+\min_{p}\frac{\lambda r_p}{1-r_p}$. Similarly, we can also derive the other side (lower bound) by using $C_1=1+\max_{p}\frac{\lambda r_p}{1-r_p}$, which eventually gives $C_1\mathcal{L}^{\Theta}_{ss}\left(\frac{\lambda-1}{C_1}\right)$.
\end{proof}}

\end{document}